\documentclass[10pt]{amsart}
\usepackage{amssymb,amsmath,amsfonts,amscd}

\usepackage[usenames]{color}
\usepackage{xcolor}

\theoremstyle{plain}
\newtheorem{theorem}{Theorem}[section]

\newtheorem{dfi}[theorem]{Definition}

\newtheorem{proposition}[theorem]{Proposition}
\newtheorem{definition}[theorem]{Definition}

\newtheorem{corollary}[theorem]{Corollary}
\newtheorem{remark}[theorem]{Remark}
\newtheorem{lemma}[theorem]{Lemma}

\usepackage{pdfsync}

\newcommand{\rn}[1]{{\mathbb R}^{#1}}
\newcommand{\R}{\mathbb R}

\newcommand{\he}[1]{{\mathbb H}^{#1}}

\newcommand{\cov}[1]{{\bigwedge\nolimits^{#1}{\mfrak h}}}

\newcommand{\covH}[1]{{\bigwedge\nolimits^{#1}{\mfrak h}}}
\newcommand{\vetH}[1]{{\bigwedge\nolimits_{#1}{\mfrak h}}}

\newcommand{\covh}[1]{{\bigwedge\nolimits^{#1}{\mfrak h_1}}}
\newcommand{\veth}[1]{{\bigwedge\nolimits_{#1}{\mfrak h_1}}}

\newcommand{\scal}[2]{\langle {#1} , {#2}\rangle}

\newcommand{\Scal}[2]{\langle {#1} \vert {#2}\rangle}
\newcommand{\scalp}[3]{\langle {#1} , {#2}\rangle_{#3}}

\newcommand{\mc}{\mathcal }

\newcommand{\mfrak}{\mathfrak}

\newcommand{\gf}{\varphi}

\newcommand{\WO}[3]{\mathop{W}\limits^\circ{}\!^{{#1},{#2}}}

\newcommand{\N}{\mathbb N}

\usepackage[usenames]{color}

\begin{document}


\title[Sobolev-Gaffney type inequalities]
{  Sobolev-Gaffney type inequalities for differential forms on sub-Riemannian contact manifolds with bounded geometry}

\author[Annalisa Baldi, Maria Carla Tesi, Francesca Tripaldi]{
Annalisa Baldi, Maria Carla Tesi  and Francesca Tripaldi
}

\begin{abstract}

In this paper
we establish a Gaffney type inequality, in $W^{\ell,p}$-Sobolev spaces, for differential forms  on  sub-Riemannian contact manifolds  without boundary, having  {bounded geometry} (hence, in particular, we have in mind non-compact manifolds). Here $p\in]1,\infty[$ and $\ell=1,2$ depending on the order of the differential form we are considering. The proof relies on the stucture of the Rumin's complex of differential forms in contact manifolds,  on a Sobolev-Gaffney inequality proved by Baldi-Franchi in the setting of the Heisenberg groups and on some geometric properties that can be proved for sub-Riemannian contact manifolds with  {bounded geometry}.  

%

\end{abstract}
 


\maketitle

\section{Introduction}

 A well known formulation of Gaffney's inequality in $\R^n$ is the following div-curl type estimate: there exists a geometric constant $C>0$ such that for any vector field  $\vec{F}$ in $W^{1,2}(\R^n,\R^n)$  
$$
\| \vec{F}\|_{W^{1,2}}\le C \|\,\mathrm{curl}\, \vec{F}\|_{L^2}+\|\,\mathrm{div}\, \vec{F}\|_{L^2}+\|\ \vec{F}\|_{L^2}\,.
$$
Such an estimate plays a
fundamental role in mathematics, for example in classical continuum and electromagnetic field theories. In the context of differential forms (by identifying 1-forms with vector fields) the previous inequality generalizes to the following, if $\alpha\in W^{1,2}(\R^n, \bigwedge^h\R^n)$,
\begin{equation}\label{tipo}
\|\alpha\|_{W^{1,2}}\le C \|d\alpha\|_{L^2}+\|\delta\alpha\|_{L^2}+\|\ \alpha\|_{L^2}\,,
\end{equation}
where $d$ and $\delta$
denote respectively the differential and codifferential of the de Rham complex in $\R^n$. Gaffney's inequality is the key estimate for the Hodge decomposition theorem for
differential forms and makes sense also in the more general framework of manifolds. Indeed, the proof of Gaffney's inequality for differential forms goes back
to Gaffney (\cite{gaffney}) for manifolds without boundary, and to Friedrichs \cite{friedrichs} and Morrey \cite{morrey} for manifolds with boundary, where the differential forms  satisfy  additional conditions on the boundary  (see  also    \cite{schwarz} and the references therein) . The proof of the above inequality in a compact Riemannian manifold without boundary, replacing the Sobolev space $ W^{1,2}$ with a more general Sobolev space $W^{1,p}$, is due to Scott \cite{scott} and its counterpart,  for compact Riemannian manifolds with boundary, is due to Iwaniec-Scott-Stroffolini \cite{ISS}.
Further contributions presenting various generalizations of Gaffney's inequality can be found in the literature. For an exhaustive overview of such Gaffney-type inequalities in the Riemannian setting we refer for example to \cite{CDS}. 

In this paper we deal with sub-Riemannian contact manifolds, and since by the classical Darboux' theorem any contact manifold is locally 
diffeomorphic to the Heisenberg group of corresponding dimension via a contact map, it is worth mentioning the specific generalization  of the  Gaffney-Friedrichs inequality recently proved in the setting of Heisenberg group in \cite{franchi_montefalcone_serra}. 

	\bigskip
	
	By replacing the exterior differential $d$  by a suitable differential $d_c^M$ which acts on differential forms ``adapted'' to the contact geometry (the so called Rumin complex, see subsection below),
	we shall prove a ${W^{\ell,p}}$-Gaffney-type inequality for a complete, non-compact sub-Riemannian contact manifold $M$ without boundary (here $\ell=1$ or $\ell=2$ depending on the degree of the differential form we are dealing with). We shall also assume that $M$ has bounded geometry, which means that, roughly speaking, there exist uniform bounds on the geometric invariants of the manifold (in the paper we adopt the definition of sub-Riemannian contact manifold with bounded geoemetry given in \cite{BFP2}). Manifolds with bounded geometry generalize the concept of compact manifolds and covering of compact manifolds. Examples of manifolds with bounded geometry are Lie groups or, more generally, homogeneous spaces.
	
	\medskip

\subsection{Sub-Riemannian contact manifolds.}
A contact manifold is given by the couple $(M,H)$ where $M$ is a smooth odd-dimensional (connected) manifold of dimension $2n+1$ and $H$ is the so called contact structure on $M$, i.e. $H$ is a smooth distribution of hyperplanes which is maximally non-integrable: given $\theta^M$ a smooth 1-form defined on $M$ such that $H=\mathrm{ker}(\theta^M)$, then $d\theta^M$ restricts to a non-degenerate 2-form on $H$.   Roughly speaking, to be maximally
non-integrable means that the contact subbundle $H$
is as far as possible from being integrable. Indeed, in general,  for a subbundle defined by a 1-form $\eta$
to be integrable it is necessary and sufficient that $\eta\wedge (d\eta)^{n}\equiv 0$ (see  \cite{mcduff_salamon}, Section 3.4). 
Therefore, a measure on a contact manifold $(M,H)$ can be defined through the non-degenerate top degree form $\theta^M\wedge(d\theta^M)^n$.
There exists, in addition,  a unique vector field $\xi^M$ transverse to $\ker \theta^M$ (the so-called Reeb vector field) such that $\theta^M(\xi^M)=1$ and $\mathcal L_{\xi^M}=0$.
 
Among stratified nilpotent Lie groups, the Heisenberg group  is the simplest example of
a group endowed with a contact structure. We recall that
the Heisenberg group $\he n$  is the Lie group with  stratified nilpotent Lie algebra $\mathfrak h$
of step 2
\begin{eqnarray*}
\mathfrak h = \mathrm{span}\,\{X_1,\dots,X_n,Y_1,\dots,Y_n\}\oplus\mathrm{span}\,\{T\}:= \mathfrak h_1\oplus 
\mathfrak h_2,
\end{eqnarray*}
where the only nontrivial commutation rules are $[X_j,Y_j]=T$, $j=1,\dots,n$. 
We denote by
$\theta^{\mathbb H }$ the 1-form on $\he n$ such that $\ker\theta^{\mathbb H} = \exp(\mathfrak h_1)$ and $\theta^{\mathbb H}(T)=1$.
We recall also that
 $\he n$ can be identified with $\rn{2n+1}$ through the exponential map. The stratification
 of the algebra induces a family of dilations $\delta_\lambda$ in the group via the exponential map,
\begin{equation}\label{uno}
\delta_\lambda=\lambda\textrm{ on }\mathfrak{h}_1,\quad \delta_\lambda=\lambda^2 \textrm{ on }\mathfrak{h}_2,
\end{equation}
which are analogous to the Euclidean homotheties.

 As already stressed above,  the Heisenberg groups are the local model of all
 contact manifolds (see \cite{mcduff_salamon}, p. 112).

Following Rumin, (see \cite{rumin_jdg} p.288) we can assume that there is a metric $g^M$ which is globally adapted to the symplectic form $d\theta^M$. Indeed,  there exists
  an endomorphisme $J$ of $H$ such that $J^2=-Id$,  $
d\theta^M(X,JY)=-d\theta^M(JX,Y)
$
for all $X, Y\in H$, 
and $d\theta^M(X,JX)>0$ for all $X\in H\setminus\{0\}$. Then, if $X,Y\in H$ we define  $g_H(X,Y):=d\theta^M(X,JY)$. Finally we extend $J$ to $TM$ by setting $J\xi^M=\xi^M$ and setting $g^M(X,Y)=\theta^M(X)\theta^M(Y)+d\theta^M(X,JY)$ for all $X,Y\in TM$.

The couple $(M,H)$ equipped with the Riemannian metric $g_H$  is called a \emph{sub-Riemannian} contact manifold and in the sequel it will be denoted by $(M,H,g^M)$, where $g^M$ is obtained as above.  In any sub-Riemannian contact manifold $(M,H,g^M)$, we can define a sub-Riemannian
distance $d_M$ (see e.g., [43]) inducing on $M$ the same topology of $M$ as a manifold. In
particular, Heisenberg groups $\he n$ can be viewed as sub-Riemannian contact manifolds. If we choose on the contact sub-bundle of $\he n$ a left-invariant metric, it turns out that
the associated sub-Riemannian metric is also left-invariant. It is customary to call this
distance in $\he n$ a Carnot-Carath\'eodory distance (notice also 
that all left-invariant sub-Riemannian metrics on Heisenberg groups are bi-Lipschitz equivalent).

A natural setting when dealing with differential forms  in Heisenberg groups  is provided by 
 Rumin's complex $(E_0^\bullet,d_c^\mathbb H)$ of differential forms in $\he n$
(see, e.g., \cite{rumin_jdg}),
since de Rham's complex $(\Omega^\bullet,d)$ in $\rn{2n+1}$, endowed with the usual exterior
differential $d$, does not fit the very structure of the group $\he n$.  Indeed, 
differential forms on $\mathfrak{h}$  split into 2 eigenspaces under $\delta_\lambda$ (see \eqref{uno}), therefore de Rham's complex lacks scale invariance under the  anisotropic dilations $\delta_\lambda$, 
 basically since it mixes derivatives along all the layers
of the stratification. 
Rumin's substitute for de Rham's exterior
differential is a linear differential operator $d_c^\mathbb H$ from sections of $E_0^h$
 to sections of $E_0^{h+1}$ ($0\le h\le 2n+1$) such
that $(d_c^\mathbb H)^2=0$.

We note explicitly 
that Rumin's differential $d_c^\mathbb H$ may be a left- invariant differential operator of order higher  than $1$.

 Rumin's construction of the intrinsic complex makes sense for arbitrary contact manifolds $(M,H)$ (see \cite{rumin_jdg}). 
 The main features of Rumin's complex defined on $M$ are the sameas those  already stated in $\he n$. If we denote by $E_0^\bullet=\oplus_{h=0}^{2n+1} E_0^h$ endowed with the exterior differential $d_c^M$, we have:
\begin{itemize}
\item[i)] $(d_c^M)^2=0$;
\item[ii)] the complex $(E_0^\bullet,d_c^M)$  is homotopically equivalent to de Rham's complex
$ (\Omega^\bullet,d)$; 
\item[iii)]  $d_c^M: E_0^h\to E_0^{h+1}$ is a homogeneous differential operator in the 
horizontal derivatives
of order 1 if $h\neq n$, whereas $d_c^M: E_0^n\to E_0^{n+1}$ is a homogeneous differential operator in the 
horizontal derivatives
of order 2.
\end{itemize}

 Given $(M,H,g^M)$, the scalar
product on $H$ determines  a norm on the
line bundle $T M/H$. Therefore for any $h$,the  $E_0^h$
are endowed with a scalar product. Using $\theta^M\wedge d\theta^M$ as
a volume form, one gets $L^p$-norms on spaces of smooth Rumin's differential forms on $M$ {(see Remark \ref{norme_lp})}.

We denote by $\ast$ the Hodge duality associated with the inner product in
$E_0^\bullet$ and the
volume form (see also Section \ref{conticodiff})and by $\delta_c^M$ the formal adjoint in $L^2(M, E_0^\bullet)$  of the operator $d_c^M$; we have $\delta^M_c=(-1)^h\ast d_c^M \ast $ on $E_0^h$ (see \cite{rumin_jdg}).

\subsection{Bounded geometry} 

In the sequel we will assume that $M$ is a non-compact manifold with bounded geometry, without boundary. But, of course, our approach covers also the case of a compact manifold without boundary.

\begin{dfi}\label{contact}
Let  $k\in\N$. Let $B(e,1)$ denote the sub-Riemannian unit  ball in $\he n$. We say that a complete smooth sub-Riemannian contact manifold $(M,H,g^M)$ has $C^k$-\emph{bounded geometry} if there exist two constants $r>0$ and $C_M=C(M)>0$ such that, for every $x\in M$, if we denote by $B(x,r)$ the sub-Riemannian ball for $(M,H,g^M)$ centered at $x$ and of radius $r$,
there exists a contactomorphism (i.e. a diffeomorphism preserving the contact forms) $\phi_x : B(e,1)\to M$ that satisfies
\begin{enumerate}
  \item $B(x,r)\subset\phi_x(B(e,1))$.
  \item $\phi_x$ is $C_M$-bi-Lipschitz.
  \item Coordinate changes $\phi_y^{-1}\circ\phi_x$ and their  derivatives up to order $k$ with respect to unit left-invariant horizontal vector fields are bounded by $C_M$.
\end{enumerate}
\end{dfi}

Examples of non compact sub-Riemannian  contact manifolds with  bounded geometry  are given  in \cite{BFP2} (see in particular Remark 4.10 therein and Section 7).

\bigskip

The formulation of a Gaffney-type inequality in this setting requires different statements, depending on the degree of the differential forms we are considering. Indeed, the fact that the operator $d_c^M$ has order 1 or 2, depending on the degree of the form on which it acts,  will be a major issue in the proofs of our results and will change the
class of Sobolev spaces in our inequalities.

To introduce the notion of Sobolev space in $M$ we will make use of the analogous notion  in Heisenberg groups $\he n$. This notion is associated with the stratification of its algebra, and nowadays is quite 
 classical: we refer, e.g., to Section 4 of \cite{folland} and Section \ref{sobolev kernels} below. Roughly speaking, 
if $\ell\in\mathbb N$ and $p\ge 1$, the Sobolev space $W^{\ell,p}(\he n)$  (the so-called Folland-Stein Sobolev
space) can be defined as follows. Fix an orthonormal basis $\{W_i,\, i=1,\dots,2n+1\} $ of $\mathfrak h$ such
that $W_i\in \mathfrak h_1$ for $i=1,\dots,2n$ and $W_{2n+1}\in\mathfrak h_2$. 
We call {\sl homogeneous order} of a monomial in $\{W_i\} $ its degree of homogeneity with respect
to the group dilations $\delta_\lambda$, $\lambda>0$.
We say that a differential form on $\he n$ belongs to $W^{\ell,p}(\he n)$ if all its derivatives of homogeneous order $\le \ell$ along $\{W_i\} $
belong to $L^p(\he n)$. 
Using $C^k$-bounded charts, this notion extends to $C^k$-bounded geometry sub-Riemannian contact manifolds $M$ (we refer to Section \ref{function spaces} for precise definitions).

Our main result reads as follows.

\begin{theorem}\label{gaffney-M-intro} Let $(M,H,g^M)$ be a complete, smooth contact manifold   with $C^k$-bounded geometry ($k\ge 3$),  without boundary. We have: 
\begin{itemize}

\item[i)] if $1\le h\le 2n$,  and $1<p<\infty$, then there
exists $C>0$ such that for all $\alpha \in  W^{1,p}(M, E_0^h)$ with $h\neq n, n+1$
\begin{equation*}\label{p>1-Mintro}\begin{split}
\| \alpha \|_{W^{1,p}(M, E_0^h)} &
\le C\big( \| d_c^M\alpha \|_{L^{p}(M, E_0^{h+1})} +
\| \delta_c^M \alpha \|_{L^{p}(M, E_0^{h-1})}
\\& \hphantom{xxxxx}+ \|\alpha\|_{L^{p}(M, E_0^h)}
\big).
\end{split}\end{equation*}
\end{itemize}

Moreover

\begin{itemize}
\item[ii)]  if $h= n$,  and $1<p<\infty$, then there
exists $C>0$ such that for all $\alpha \in  W^{2,p}(M, E_0^n)$
\begin{equation*}\label{treno1intro}\begin{split}
\| \alpha \|_{W^{2,p}(M, E_0^n)} &
\le C\big( \| d_c^M\alpha \|_{L^{p}(M, E_0^{n+1})} +
\| d_c^M\delta_c^M \alpha \|_{L^{p}(M, E_0^{n})}
\\&
 \hphantom{xxxxx} +  \|\alpha\|_{L^p(M, E_0^n)}\big).
\end{split}\end{equation*}
\end{itemize}

On the other hand, \begin{itemize}
\item[iii)]  if $h=  n+1$,   and $1<p<\infty$, then there
exists $C>0$ such that for all  $\alpha \in  W^{2,p}(M, E_0^n)$
\begin{equation*}\label{p>1:eq2Mintro}\begin{split}
\| \alpha \|_{W^{2,p}(M, E_0^{n+1})} &
\le C\big( \| \delta_c^Md_c^M\alpha \|_{L^{p}(M, E_0^{n+1})} +
\| \delta_c^M \alpha \|_{L^{p}(M, E_0^{n})}
\\&
 \hphantom{xxxxx} +  \|\alpha\|_{L^p(M, E_0^{n+1})}\big).
\end{split}\end{equation*}
\end{itemize}
 \end{theorem}


\bigskip

The paper is organized as follows. Section 2 contains some definitions and properties of Heisenberg groups and a very short review of Rumin's complex.  Section 3  contains a precise definition of manifold with bounded geometry and the  definition of its associated  Sobolev spaces. Section 3  also contains a fine result concerning the existence of a map that associates  an orthonormal symplectic basis of $\ker\theta^M$ with the  canonical orthonormal symplectic basis of $\mathrm{Ker \,} \theta^{\mathbb H}_e$, where $\theta^{\mathbb H}_e$ denotes the form $\theta^{\mathbb H}$ evaluted at the point $e\in \he n$. It should be noted that throughout the paper we will use both notations $\alpha_p$ and $\alpha(p)$ to indicate when the form $\alpha$ is evalueted at the point $p$.   In our setting, if  $\psi :\he n\to M$ is a contact diffeomorphism,  we have  $\psi^\sharp d_c^M=d_c^\mathbb H\psi^\sharp$ (where $\psi^\sharp$ is the pullback of the map $\psi$).  This is not the case when we replace the differential with the codifferential; in Section 4 we show that if $M$ has bounded geometry  we can locally control the difference $\psi^\sharp \delta_c^M-\delta_c^\mathbb H\psi^\sharp$ with constants depending only on the geometry of $M$. In Section 5 we prove at first a local Gaffney inequality, and then we pass from the local case to the global one using an   ad hoc  covering for M with balls of  suitably small radius obtained in Section 3.

We now list  a few difficulties we
have to deal with in this note. To prove the Gaffney inequality on $M$, we pass from
the corresponding result for the ``flat'' model $\he n$  proved in \cite{BF7}.  A delicate point is to show that we can replace the contactomorphism  $\phi_x$ appearing in Definition \ref{contact} with another  contactomorphism $\psi_x$ which ``sends'' an orthonormal symplectic basis of $\ker\theta^M_x$ into the  canonical orthonormal symplectic basis of $\ker \theta^{\mathbb H}_e$  and that still depends only on the bounded geometric constants and not on the point $x$.  This is accomplished in Theorem \ref{delirio} and it is a convenient step if we want to write in  local coordinates the difference $\psi^\sharp \delta_c^M-\delta_c^\mathbb H\psi^\sharp$. Another subtle issue comes from the fact that the order of  the differentials $d_c^M$ and $d_c^\mathbb H$ can be one or two depending on the degree of the form and this problem is reflected in several proofs (for example in Section 4, when we estimate $\psi^\sharp \delta_c^M-\delta_c^\mathbb H\psi^\sharp$, or in Theorem \ref{gaffney-M} when we have to estimate the commutators between differentials  and functions).

\section{Basic properties of Rumin's complex $(E_0^\bullet,d_c)$ on  Heisenberg groups and on general contact manifolds}
\label{heisenberg}

\subsection{Heisenberg groups}
 We denote by  $\he n$  the { $2n+1$}-dimensional Heisenberg
group, identified with $\rn {2n+1}$ through exponential
coordinates. A point $p\in \he n$ is denoted by
 by 
$p=(x,y,t)$, with both $x,y\in\rn{n}$
and $t\in\R$.
   If we consider two points  in $ \he n$, $p=(x,y,t)$ and $q=(\tilde x, \tilde y, \tilde t)$
,   the (non commutative) group operation is denoted by $p\cdot q:=(x+\tilde x, y+\tilde y, t+\tilde t + \frac12 \sum_{j=1}^n(x_j \tilde y_{j}- y_{j} \tilde x_{j}))$.
In this system of coordinates,
the  unit element of $\he n$, that will be denote by $e$, is the zero of the vector space $\R^{2n+1}$, and $p^{-1} = -p$.

For
any $q\in\he n$, the {\it (left) translation} $\tau_q:\he n\to\he n$ is defined
as $$ p\mapsto\tau_q p:=q\cdot p. $$
The Lebesgue measure in $\mathbb R^{2n+1}$ 
is a Haar measure in $\he n$.

For a general review on Heisenberg groups and their properties, we
refer to \cite{Stein}, \cite{GromovCC} and to \cite{VarSalCou}.
We limit ourselves to fix some notations, following e.g. \cite{BFP2}.

First we notice that Heisenberg groups are smooth manifolds (and therefore
are Lie groups). In particular, the pullback of differential forms is well 
defined as follows (see, e.g. \cite{GHL}, Proposition 1.106);
\begin{definition} If $\; \mc U,\mc V$ are open subsets of $\he n$, and $f: \mc U\to
\mc V$ is a
diffeomorphism, then for any  differential form  $\alpha$ of degree $h$,
 we denote by $f^\sharp \alpha$ the pullback form in $\mc U$ defined by
$$
(f^\sharp \alpha)(p) (v_1,\dots,v_h):=  \alpha(f(p)) (df(p)v_1,\dots,df(p)v_h)
$$
for any $h$-tuple $(v_1,\dots,v_h)$ of tangent vectors at $p$.

\end{definition}

The Heisenberg group $\he n$ can be endowed with the homogeneous
norm (Cygan-Kor\'anyi norm), 
$$\varrho (p)=\big((|x|^2+|y|^2)^2+ 16t^2\big)^{1/4}$$
and we define the gauge distance (a true distance, see
 \cite{Stein}, p.\,638)
as
\begin{equation}\label{def_distance}
d(p,q):=\varrho ({p^{-1}\cdot q}).
\end{equation}
The metric $d$  behaves wel lwith respect to left-translations,
that is
$$d(\tau_q p,\tau_q p') = d(p,p')\,, $$
for all $q, p,p'\in\he n$.
Finally, the balls for the metric $d$ are the so-called Kor\'anyi balls
\begin{equation}\label{koranyi}
B(p,r):=\{q \in  \he n; \; d(p,q)< r\}.
\end{equation}

    We denote by  $\mfrak h$
 the Lie algebra of the left-invariant vector fields of $\he n$. The standard basis of $\mfrak
h$ is given, for $i=1,\dots,n$,  by
\begin{equation*}
X_i := \partial_{x_i}-\frac12 y_i \partial_{t},\quad Y_i :=
\partial_{y_i}+\frac12 x_i \partial_{t},\quad T :=
\partial_{t}.
\end{equation*}
The only non trivial bracket  relations are $
[X_{j},Y_{j}] = T $, for $j=1,\dots,n.$ 
The {\it horizontal subspace}  $\mfrak h_1$ is the subspace of
$\mfrak h$ spanned by $X_1,\dots,X_n$ and $Y_1,\dots,Y_n$.
Coherently, from now on, we refer to $X_1,\dots,X_n,Y_1,\dots,Y_n$
(identified with first order differential operators) as
the {\it horizontal derivatives}. Denoting  by $\mfrak h_2$ the linear span of $T$, the $2$-step
stratification of $\mfrak h$ is expressed by
\begin{equation*}
\mfrak h=\mfrak h_1\oplus \mfrak h_2.
\end{equation*}

\bigskip

The stratification of the Lie algebra $\mfrak h$ induces a family of non-isotropic dilations
$\delta_\lambda$, $\lambda>0$ in $\he n$ so that for any $p=(x,y,t)\in \he n$ 
\begin{equation}\label{dilatazione}
	 \delta_\lambda(p)=(\lambda x, \lambda y,\lambda^2 t)\,.
\end{equation}


Notice that the gauge norm  is positively $\delta_\lambda$-homogenous
(i.e. $
 d(\delta_\lambda(p),\lambda(p')) = \lambda d(p, p')
$ for all $q, p,p'\in\he n$ and $\lambda > 0$)
so that the Lebesgue measure of the ball $B(x,r)$ is $r^{2n+2}$ up to a geometric constant
(the Lebesgue measure of $B(e,1)$).
Thus, the {\sl homogeneous dimension}  of $\he n$
with respect to $\delta_\lambda$, $\lambda>0$, equals 
$$
Q:=2n+2.
$$
It is well known that the topological dimension of $\he n$ is $2n+1$,
since as a smooth manifold it coincides with $\R^{2n+1}$, whereas
the Hausdorff dimension of $(\he n,d)$ is $Q$.

The vector space $ \mfrak h$  can be
endowed with an inner product, indicated by
$\scalp{\cdot}{\cdot}{} $,  making
    $X_1,\dots,X_n$,  $Y_1,\dots,Y_n$ and $ T$ orthonormal.
    
Throughout this paper, we write also
\begin{equation}\label{campi W}
W_i^{\mathbb H}:=X_i, \quad W_{i+n}^{\mathbb H}:= Y_i, \quad W_{2n+1}^{\mathbb H}:= T, \quad \text
{for }i =1, \cdots, n.
\end{equation}

As in \cite{folland_stein},
we also adopt the following multi-index notation for higher-order derivatives. If $I =
(i_1,\dots,i_{2n+1})$ is a multi--index, we set  
\begin{equation}\label{WI}
W^{\mathbb H, I}=(W_1^\mathbb H)^{i_1}\cdots
(W_{2n}^\mathbb H)^{i_{2n}}\;(W_{2n+1}^{\mathbb H})^{i_{2n+1}}.
\end{equation}
By the Poincar\'e--Birkhoff--Witt theorem, the differential operators $W^{\mathbb H,I}$ form a basis for the algebra of left invariant
differential operators in $\he n$. 
Furthermore, we set 
$$
|I|:=i_1+\cdots +i_{2n}+i_{2n+1}$$
 the order of the differential operator
$W^{\mathbb H,I}$, and   
$$d(I):=i_1+\cdots +i_{2n}+2i_{2n+1}$$
 its degree of homogeneity
with respect to group dilations.

\subsubsection{Sobolev spaces in $\he n$}\label{sobolev kernels}

Let  $U\subset \he n$ be an open set. We shall use the following classical notations:
$\mc E(U)$ is the space of all smooth function on $U$,
and $\mc D(U)$  is the space of all compactly supported smooth functions on $U$,
endowed with the standard topologies (see e.g. \cite{treves}).

We recall  the notion of 
(integer order) Folland-Stein Sobolev space (for a general presentation, see e.g. \cite{folland} and \cite{folland_stein}).

\begin{definition}\label{integer spaces} If $U\subset \he n$ is an open set, $1\le p \le\infty$
and $k\in\mathbb N$, then
the space $W^{k,p}(U)$
is the space of all $u\in L^p(U)$ such that, with the notation of \eqref{WI},
$$
W^{\mathbb H, I}u\in L^p(U)\quad\mbox{for all multi-indices $I$ with } d(I) \le k,
$$
endowed with the natural norm  
$$\|\,u\|_{W^{k,p}(U)}
:= \sum_{d(I) \le k} \|W^{\mathbb H, I} u\|_{L^p(U)}.
$$

\end{definition}

{  Folland-Stein Sobolev spaces enjoy the following properties akin to those
of the usual Euclidean Sobolev spaces (see \cite{folland}, and, e.g. \cite{FSSC_houston}).}

\begin{theorem}\label{denso in h} If $U\subset \he n$, $1\le p \le  \infty$, and $k\in\mathbb N$, then
\begin{itemize}
\item[i)] $ W^{k,p}(U)$ is a Banach space.
\end{itemize}
In addition, if $p<\infty$,
\begin{itemize}
\item[ii)] $ W^{k,p}(U)\cap \mc E(U)$ is dense in $ W^{k,p}(U)$;
\item[iii)] if $U=\he n$, then $\mc D(\he n)$ is dense in $ W^{k,p}(\he n)$;
\item[iv)]  if $1<p<\infty$, then $W^{k,p}(U)$ is reflexive.
\end{itemize}

\end{theorem}

\begin{definition} { If $U\subset \he n$ is open} and if
$1\le p<\infty$,
we denote  by $\WO{k}{p}{U}$
the completion of $\mc D(U)$ in $W^{k,p}(U)$.
\end{definition}
{  \begin{remark}
If $U\subset \he n$ is bounded, by (iterated) Poincar\'e inequality (see e.g. \cite{jerison}), it follows that the norms
\begin{equation*}
\|u\|_{W^{k,p}(U)} \quad\mbox{and}\quad
\sum_{d(I)=k}\| W^{\mathbb H ,I} u\|_{ L^p(U)}
\end{equation*}
are equivalent on $\WO{k}{p}{U}$ when $1\le p<\infty$.
\end{remark}
}

We recall the following Ehrling's inequality.
%
\begin{remark}\label{magenes} 
Let $B(e,1)$ be the Kor\'anyi ball in $\he n$. By \cite{MoMo}  $B(e,1)$ is a $(\epsilon,\delta)$ domain, hence we have  the compact  embedding $W^{1,p}(B(e,1))\hookrightarrow L^{p}(B(e,1))$ by  e.g.  Theorem 1.27 in \cite{GN} (see also \cite{lu_acta_sinica}, \cite{hajlasz_koskela_memoires}). In conclusion, $W^{2,p}(B(e,1))\hookrightarrow W^{1,p}(B(e,1))\hookrightarrow L^{p}(B(e,1))$.  Hence, by Ehrling's inequality (see e.g. 
\cite{schwarz} Lemma 1.5.3) ,
if $v\in W^{2,p}(B(e,1))$ , {{for any $0<\varepsilon<1$ there exists a constant $c(\varepsilon)$ such that }}we have:
\begin{equation*}
	\|v\|_{W^{1,p}(B(e,1))}\le \varepsilon \|v\|_{W^{2,p}(B(e,1))}+c(\varepsilon) \|v\|_{L^{p}(B(e,1))}\,
\end{equation*}
\end{remark}
\subsection{Multilinear algebra and Rumin's complex in Heisenberg groups}

The dual space of $\mfrak h$ is denoted by $\covH 1$.  The  basis of
$\covH 1$,  dual to  the basis $\{X_1,\dots , Y_n,T\}$,  is the family of
covectors $\{dx_1,\dots, dx_{n},dy_1,\dots, dy_n,\theta^{\mathbb H}\}$ where 
$$ \theta^{\mathbb H}
:= dt - \frac12 \sum_{j=1}^n (x_jdy_j-y_jdx_j)$$ is called the {\it contact
form} in $\he n$.

We indicate again by $\scalp{\cdot}{\cdot}{} $  the
inner product in $\covH 1$  that makes $(dx_1,\dots, dy_{n},\theta^{\mathbb H}  )$ 
an orthonormal basis. Coherently with the previous notation \eqref{campi W},
we set
\begin{equation*}
\omega_i^{\mathbb H}:=dx_i, \quad \omega_{i+n}^{\mathbb H}:= dy_i, \quad \omega^{\mathbb H}_{2n+1}:= \theta^{\mathbb H}, \quad \text
{for }i =1, \cdots, n.
\end{equation*}
The volume $(2n+1)$-form $ \omega_1^{\mathbb H}\wedge\cdots\wedge \omega^{\mathbb H}_{ 2n+1}$
 will be also
written as $dV$.

We denote by 
$       \vetH 0 := \covH 0 =\R $
and, for $1\leq h \leq 2n+1$, we can define
\begin{equation*}
\begin{split}
 \vetH h& :=\mathrm {span}\{ W_{i_1}^{\mathbb H}\wedge\dots \wedge W_{i_h}^{\mathbb H}: 1\leq
i_1< \dots< i_h\leq 2n+1\},   \\
         \covH h& :=\mathrm {span}\{ \omega_{i_1}^{\mathbb H}\wedge\dots \wedge \omega^{\mathbb H}_{i_h}:
1\leq i_1< \dots< i_h\leq 2n+1\}
.
\end{split}
\end{equation*}

In the sequel we shall denote by $\Theta^h$ the basis of $ \covH h$ defined by
$$
\Theta^h:= \{ \omega_{i_1}^{\mathbb H}\wedge\dots \wedge \omega_{i_h}^{\mathbb H}:
1\leq i_1< \dots< i_h\leq 2n+1\}.
$$
The  {{ inner}} product $\scal{\cdot}{\cdot}$ on $ \covH 1$ yields naturally an {{ inner}} product 
$\scal{\cdot}{\cdot}$ on $ \covH h$
making $\Theta^h$ an orthonormal basis.

If $1\leq h\leq 2n+1$, the Hodge isomorphism
\begin{equation*}
\ast : \vetH h \longleftrightarrow \vetH{2n+1-h} \quad \text{and}
\quad \ast : \covH h \longleftrightarrow \covH{2n+1-h},
\end{equation*}
is defined by 
$$
v\wedge \ast w= \scal v w{} W_1^{\mathbb H}\wedge\cdots\wedge W^{\mathbb H}_{ 2n+1}\quad\quad \forall\  v,w\in\vetH h,$$
$$\gf \wedge \ast \psi= \scal \gf \psi{} \omega_1^{\mathbb H}\wedge\cdots\wedge \omega^{\mathbb H}_{ 2n+1}\quad\quad \forall\ v,w\in\covH h.$$

 If $v\in \vetH h$ we define its dual $v^\natural \in \covH h$ by the identity
$ \Scal {v^\natural} w := \scal v w , $ and analogously we define
$\gf^\natural\in \vetH h$ for $\gf \in \covH h$.

Throughout this paper, the elements of $\cov h$ are identified with \emph{left invariant} differential forms
of degree $h$ on $\he n$. 

\begin{definition}\label{left} A $h$-form $\alpha$ on $\he n$ is said left invariant if 
$$\tau_q^\sharp\alpha
=\alpha\qquad\mbox{for any $q\in\he n$.}
$$
Here $\tau_q^\sharp\alpha$ denotes the pull-back of $\alpha$ through the left translation $\tau_q$.
\end{definition}

The same construction given above can be performed starting from the vector
subspace $\mfrak h_1\subset \mfrak h$,
obtaining the {\it horizontal $h$-covectors} 
and {\it horizontal
$h$-vectors}
\begin{equation*}
\begin{split}
         \veth h& :=\mathrm {span}\{ W_{i_1}^{\mathbb H}\wedge\dots \wedge W_{i_h}^{\mathbb H}:
1\leq
i_1< \dots< i_h\leq 2n\}   \\
         \covh h& :=\mathrm {span}\{ \omega_{i_1}^{\mathbb H}\wedge\dots \wedge \omega_{i_h}^{\mathbb H}:
1\leq i_1< \dots< i_h\leq 2n\}.
\end{split}
\end{equation*}
and $$
\Theta^h_0 := \Theta^h \cap  \covh h
$$ 
provides an orthonormal
basis of $ \covh h$.

The {\it symplectic 2-form}  $d\theta^{\mathbb H} \in \covh 2$ is 
$$d \theta^{\mathbb H} =
-\sum _{i=1}^n \omega_i^{\mathbb H} \wedge \omega_{i+n}^{\mathbb H}.$$

\medskip

Keeping in mind that the Lie algebra $\mathfrak h$ can be identified with the
tangent space to $\he n$ at $x=e$ (see, e.g. \cite{GHL}, Proposition 1.72), 
starting from $\cov h$ we can define by left translation  a fiber bundle
over $\he n$  that we can still denote by $\cov h$. We can think of $h$-forms as sections of 
$\cov h$. We denote by $\Omega^h$ the
vector space of all smooth $h$-forms.

In addition, the  symplectic 2-form   
$-d \theta^{\mathbb H} $
induces on $\mathfrak h_1$  a symplectic structure.
Notice that $\{W_1^{\mathbb H}, \cdots , W_{2n}^{\mathbb H} \}$  is a symplectic basis of $\ker \theta^{\mathbb H}$.

\subsubsection{Rumin's complex}

{ Unfortunately, when dealing 
with differential forms in $\he n$, 
the de Rham complex lacks scale invariance under anisotropic dilations (see \eqref{uno}). 
Thus, a substitute for de Rham's complex, that recovers scale invariance under $\delta_\lambda$ has been defined 
by M. Rumin, \cite{rumin_jdg}. 

Here we give only a short introduction to Rumin's complex. For a more detailed presentation we
refer to Rumin's papers \cite{rumin_grenoble} or we can follow  the presentation of \cite{BFTT}. 

Let $L: \cov h \to \cov{h+2}$ the Lefschetz operator defined by
\begin{equation}\label{lefs}
L\, \xi = d\theta^{\mathbb H}\wedge\xi.
\end{equation}
Then the spaces $E_0^\bullet\subset \cov \bullet$  can be ``defined'' explicitly as follows:

\begin{theorem}[see \cite{rumin_jdg}, \cite{rumin_gafa}] \label{rumin in H} We have:
\begin{itemize}
\item[i)] $E_0^1= \covh{1}$;
\item[ii)]  if $2\le h\le n$, then $E_0^h= \covh{h}\cap \big(\covh{h-2}\wedge d\theta^{\mathbb H}\big)^\perp$
 (i.e. $E_0^h$ is the space of the so-called \emph{primitive covectors} of $\covh h$);
\item[iii)]  if $n< h\le 2n+1$, then $E_0^h = \{\alpha = \beta\wedge\theta^{\mathbb H}, \; \beta\in \covh{h-1},
\; \beta\wedge d\theta^{\mathbb H} =0\} = \theta^{\mathbb H}\wedge\ker L$;
\item[iv)]  if $1<h\le n$, then $\dim E_0^h = \binom{2n}{h} - \binom{2n}{h-2}$;
\item[v)]  if $\ast$ denotes the Hodge duality associated with the {{ inner}} product in $\cov{\bullet}$
and the volume form, then $\ast E_0^h = E_0^{2n+1-h}$.
\end{itemize} 
\end{theorem}


\medskip

For $h=0,\ldots,2n+1$, the space of Rumin $h$-forms,  is the space of smooth 
sections of a left-invariant subbundle of $\bigwedge^h\mathfrak h$, that we still denote by $E_0^{h}$. 
Hence it inherits the inner product and the norm of $\bigwedge^h\mathfrak h$.

The core of Rumin's theory consists in the construction of a suitable ``exterior differential''
$d_c^{\mathbb H}: E_0^h \to E_0^{h+1}$ making $\mc E_0:= (E_0^\bullet, d_c^{\mathbb H})$ a complex homotopic
to the de Rham complex (i.e. $d_c^{\mathbb H}\circ d_c^{\mathbb H}=0$). 

More precisely,
 the exterior differential  $d_c^{\mathbb H}:E_0^{h}\to E_0^{h+1}$  is a left-invariant, 
homogeneous operator with respect to group dilations. It is a first order homogeneous operator
in the horizontal derivatives in degree $\neq n$, whereas {\emph it is a second
order homogeneous horizontal operator in degree $n$}. 
There exists a left invariant orthonormal basis 
of $E_0^h$, such a basis is explicitly constructed by {  induction} in \cite{BBF}. 
Explicit computations  of the classes $E_0^h$ and the differential
$d_c^{\mathbb H}: E_0^h \to E_0^{h+1}$ in $\he 1$ and $\he 2$ are given for example in \cite{BF7} (see Example 3.11 and Example 3.12 therein).

The next remarkable property of Rumin's complex is its invariance under contact transformations. 
In particular,

\begin{proposition}\label{pull} If we write a form $\alpha = \sum_I \alpha_I \xi^\mathbb H_I$ in coordinates with respect to a left invariant basis $\{\xi^\mathbb H_I\}_I$ of $E_0^h$
 we have:
\begin{equation}\label{pull trasl}
\tau_q^\sharp \alpha = \sum_I ( \alpha_I \circ \tau_q)\xi^\mathbb H_I
\end{equation}
for all $q\in \he n$.
In addition, for $t>0$,
\begin{equation}\label{pull dil 1}
\delta_t^\sharp \alpha =  t^h \sum_I ( \alpha_I \circ \delta_t)\xi^\mathbb H_I\qquad\mbox{if $1\le h \le n$}
\end{equation}
and 
\begin{equation}\label{pull dil 2}
\delta_t^\sharp \alpha =  t^{h+1} \sum_I ( \alpha_I \circ\delta_t)\xi^\mathbb H_I\qquad\mbox{if $n+1\le h \le 2n+1$}\,.
\end{equation}
\end{proposition}

Let us list a list of notations for vector-valued function spaces (for the scalar case,
we refer to Section \ref{sobolev kernels}).

\begin{definition} \label{dual spaces forms-no} If $U\subset \he n$ is an open set, $0\le h\le 2n+1$, $1\le p\le \infty$ and $m\ge 0$,
we denote by $L^{p}(U,\cov{h})$, $\mc E (U,\cov{h})$,  $\mc D (U,\cov{h})$, 
$W^{m,p}(U,\cov{h})$ (by $\WO{m}{p}{U,\cov{h}}$)
the space of all sections of $\cov{h}$ such that their
components with respect to a given left invariant frame  belong to the corresponding scalar spaces.

The spaces $L^{p}(U,E_0^h)$, $\mc E (U,E_0^h)$,  $\mc D (U, E_0^h)$,  
$W^{m,p}(U,E_0^{h})$ and $\WO{m}{p}{U,E_0^h}$ are defined in the same way.

Clearly, all these definitions
are independent of the choice of frame.
\end{definition}

In addition, Sobolev spaces of differential forms are invariant with
respect to the pullback operator associated with contact diffeomorphisms (see \cite{BFP2}, Lemma 4.8).

\begin{proposition} Denote by $\delta_c^{\mathbb H} $
the formal adjoint of $d_c^{\mathbb H}$ in $L^2(\he n,E_0^h)$. 
Then $\delta_c^{\mathbb H} = (-1)^h \ast d_c^{\mathbb H}\ast$ on $E_0^h$.
\end{proposition}
We remind that $\delta_c^{\mathbb H}$ can be written in coordinates as a left-invariant homogeneous differential operator in the horizontal variables, of order 1 if $h\neq n+1$ and of order
2 if $h=n+1$ (see Example 3.11 and Example 3.12 \cite{BF7} for explicit expressions of the codifferential).

When $d_c^{\mathbb H}$ is second order (i.e. when $d_c^\mathbb H$ acts on forms of degree $n$), the complex $(E_0^\bullet,d_c^{\mathbb H})$ stops behaving like a differential module. This is the source of many complications. In particular, the classical Leibniz formula for the de Rham complex $d(\alpha\wedge\beta)=d\alpha\wedge\beta\pm \alpha\wedge d\beta$ in general fails to hold (see \cite{BFT2}-Proposition A.7). This causes several technical difficulties when we want to localize our estimates by means of cut-off functions.
}

If $\zeta$ is a smooth real function and $\alpha\in L^1_{\mathrm{loc}}(\he n, E_0^h)$ we write $d_c^\mathbb H(\zeta\alpha)=\zeta d_c^\mathbb H\alpha+[d_c^{\mathbb H},\zeta]\alpha$. The proof of the following Leibniz' formula can be found in  \cite{BFP3}, Lemma 4.1 and basically is due to the fact that  the exterior differential $d_c^{\mathbb H}$ on $E_0^h$ can be written in coordinates as a left-invariant homogeneous differential operator in the horizontal variables, of order 1 if $h\neq n$ and of order
2 if $h=n$ (and the codifferential $\delta_c^{\mathbb H}$ can be written in coordinates as a left-invariant homogeneous differential operator in the horizontal variables, of order 1 if $h\neq n+1$ and of order
2 if $h=n+1$).

\begin{lemma} \label{leibniz} If $\zeta$ is a smooth real function, then the following formulae hold:
\begin{itemize}
\item[i)] if $h\neq n$, then on $E_0^h$ we have
$$
[d_c^{\mathbb H},\zeta] = P_0^h(W\zeta),
$$
where $P_0^h(W\zeta): E_0^h \to E_0^{h+1}$ is a linear homogeneous differential operator of order zero with coefficients depending
only on the horizontal derivatives of $\zeta$. If $h\neq n+1$, an analogous statement holds if we replace
$d_c^\mathbb H$ { in degree $h$ with $\delta_c^\mathbb H$ in degree $h+1$};
\item[ii)] if $h= n$, then on $E_0^n$ we have
$$
[d_c^{\mathbb H},\zeta] = P_1^n(W\zeta) + P_0^n(W^2\zeta) ,
$$
where $P_1^n(W\zeta):E_0^n \to E_0^{n+1}$ is a linear homogeneous differential operator of order 1 (and therefore
horizontal) with coefficients depending
only on the horizontal derivatives of $\zeta$, and where $P_0^h(W^2\zeta): E_0^n \to E_0^{n+1}$ is a linear homogeneous differential operator in
the horizontal derivatives of order 0
 with coefficients depending
only on second order horizontal derivatives of $\zeta$. If $h = n+1$, an analogous statement holds if we replace
$d_c^\mathbb H$ { in degree $n$ with $\delta_c^\mathbb H$ in degree $n+1$}.

\item[iii)] if $h\neq n+1$, then
$$
[d_c^{\mathbb H}\delta_c^{\mathbb H},\zeta] = P_1^h(W\zeta) + P_0^h(W^2\zeta) ,
$$
where $P_1^h(W\zeta):E_0^h \to E_0^{h}$ is a linear homogeneous differential operator of order 1 (and therefore
horizontal) with coefficients depending
only on the horizontal derivatives of $\zeta$, and where $P_0^h(W^2\zeta): E_0^h \to E_0^{h}$ is a linear homogeneous differential operator in
the horizontal derivatives of order 0
 with coefficients depending
only on second order horizontal derivatives of $\zeta$.
\item[iv)] if $h\neq n$, then
$$
[\delta_c^{\mathbb H} d_c^{\mathbb H},\zeta] = P_1^h(W\zeta) + P_0^h(W^2\zeta) ,
$$
where $P_1^h(W\zeta):E_0^h \to E_0^{h}$ is a linear homogeneous differential operator of order 1 (and therefore
horizontal) with coefficients depending
only on the horizontal derivatives of $\zeta$, and where $P_0^h(W^2\zeta): E_0^h \to E_0^{h}$ is a linear homogeneous differential operator in
the horizontal derivatives of order 0
 with coefficients depending
only on second order horizontal derivatives of $\zeta$.

\end{itemize}

\end{lemma}

\begin{remark}\label{coleibniz} 
On forms of degree $h>n$, Lemma \ref{leibniz} i) takes the following simpler form. If $\alpha\in L^1_{\mathrm{loc}} (\he n, E_0^h)$ with $h>n$ and $\psi\in\mc E(\he n)$, then
\begin{equation*}\begin{split}
d_c^{\mathbb H}(\psi\alpha) = d(\psi\alpha) = d\psi\wedge\alpha + \psi d\alpha =
d_c^{\mathbb H}\psi\wedge\alpha + \psi d_c^{\mathbb H}\alpha,
\end{split}\end{equation*}
(in the sense of distributions). This follows from Theorem \ref{rumin in H}, iii), since $\alpha$ is a multiple of $\theta$. 
\end{remark}

\subsubsection{Rumin's complex in contact manifolds}
The notion of Rumin's complex  makes sense for arbitrary contact manifolds. 

Let us start with the following definition (see \cite{mcduff_salamon}, Section I-3).
\begin{definition}\label{uggioso} If $(M_1, H_1)$ and $(M_2,H_2)$ are contact manifolds with $H_i = \ker\theta^{M_i}$ (i.e. $\theta^{M_i}$ are contact forms) $i=1,2$,
$\mc U_1\subset M_1$, $\mc U_2\subset M_2$ are open sets and $f$ is a diffeomorphism from $\mc U_1$ onto $\mc U_2$, then $f$ is said a contact diffeomorphism if 
there exists a non-vanishing real function $\tau$ defined in $\mc U_1$ such that
$$
 f^\sharp \theta^{M_2}  =  \tau\theta^{M_1} \qquad\mbox{in \,$\mc U_1$.}
$$
\end{definition}
As already pointed out, by a classical Darboux' theorem, any contact manifold $(M,H)$  of dimension $2n+1$ is locally 
contact diffeomorphic to the Heisenberg group $\he n$.
 It turns out that Rumin's intrinsic complex is invariant under contactomorphisms, hence it is invariantly defined for general contact manifolds $(M,H)$.


A detailed  construction of this intrinsic complex for general contact manifolds $(M,H)$ can be found in  \cite{rumin_jdg} (see also  \cite{BFP2} Section 3.3 for details). 
Alternative contact invariant definitions of Rumin's complex can be found in \cite{Bernig_2017} and \cite{BEGN}.

In this paper we do not enter into the details of Rumin's construction,   we just recall  the basic properties enjoyed by  the complex, exactly analogous to the ones in $\he n$.
\begin{itemize}
\item[i)] $d_c^M\circ d_c^M=0$;
\item[ii)] the complex $\mc E_0:=(E_0^\bullet,d_c^M)$ is homotopically equivalent to the de Rham complex
$ (\Omega^\bullet,d)$. 
\item[iii)] $d_c^M: E_0^h\to E_0^{h+1}$ is a  differential operator 
of order 1 if $h\neq n$, whereas $d_c^M: E_0^n\to E_0^{n+1}$ is a differential operator  
of order 2.
\end{itemize}

In particular we recall the following statement proved in \cite{BFP2}, Proposition 2.14, that we shall use in the sequel and that  expresses the fact that Rumin's complex is invariant under contactomorphism. 

\begin{proposition}\label{memory} If $\phi $ is a contactomorphism  from an open set $\mc U\subset \he{n}$ to $M$, and we denote by $\mc V$
the open set  $\mc V:= \phi(\mc U)$, the pullback operator $\phi^\sharp$ satisfies:
\begin{itemize}
\item[i)] $\phi^\sharp E_0^\bullet(\mc V) = E_0^\bullet(\mc U) $;
\item[ii)] $d_c^\mathbb H\phi^\sharp = \phi^\sharp d_c^M$.
\end{itemize}

\end{proposition}

\section{Sobolev spaces on contact sub-Riemannian manifolds with bounded geometry}\label{function spaces}

We define Sobolev spaces (involving  a positive  number of derivatives) on  contact sub-Riemannian manifolds with bounded geometry.

We make  the definition  of contact manifold of bounded geometry, already given in  Definition \ref{contact}, more precise. 
\begin{dfi}\label{contact bis} Let $k$ be a positive  integer and
let $B(e,1)$ denote the unit sub-Riemannian ball in $\he n$.  We say that a sub-Riemannian contact manifold $(M,H,g^M)$ has \emph{bounded $C^k$-geometry} if there exist constants $r, C_M>0$ such that, for every $x\in M$, 
there exists a contactomorphism (i.e. a diffeomorphism preserving the contact forms) $\phi_x : B(e,1)\to M$ {that satisfies}
\begin{enumerate}
  \item $B(x,r)\subset\phi_x(B(e,1))$;
  \item $\phi_x$ is $C_M$-bi-Lipschitz, i.e.
  \begin{equation}\label{bilip}
\frac{1}{C_M}d(p,q)\le  d_M(\phi_x (p), \phi_x(q)) \le C_Md(p,q)\qquad \mbox{for all $p,q\in B(e,1)$};
\end{equation}
  \item coordinate changes $\phi_y^{-1}\circ \phi_x$ and their first $k$ derivatives with respect to unit left-invariant horizontal vector fields are bounded by $C_M$.
\end{enumerate}
\end{dfi}

We recall the following covering lemma  and definition from \cite{BFP2}. 
\begin{lemma}\label{covering}[see \cite{BFP2}, Lemma 4.11] Let $(M,H,g^M)$ be a  $C^k$-bounded geometry sub-Riemannian contact manifold,
where $k$ is a positive integer. Then there
exists $\rho>0$ (depending only on the radius $r$ of Definition \ref{contact bis}) and an at most countable covering $\{ B(x_j, \rho)\}$ of $M$ such that
\begin{itemize}
\item[i)] each ball $B(x_j, \rho)$ is contained in the image of one of the contact charts of Definition \ref{contact bis};
\item[ii)] $B(x_j, \frac15\rho)\cap B(x_i, \frac15\rho)=\emptyset$ if $i\neq j$;
\item[iii)] the covering is uniformly locally finite. Even more,  there exists $N=N(M)\in \mathbb N$ such that 
for each ball $B(x,\rho)$ 
$$
\#\{k\in \mathbb N \mbox{ such that } B(x_k,\rho)\cap B(x,\rho)\neq\emptyset\} \le N.
$$
In addition, if $B(x_k,\rho)\cap B(x,\rho)\neq\emptyset$, then $B(x_k,\rho)\subset B(x,r)$, where $B(x,r)$ has
been defined in Definition \ref{contact bis}).
\end{itemize}

\end{lemma}

We are in position to define Sobolev spaces on $M$  on bounded geometry contact sub-Riemannian manifolds.
\begin{definition}\label{carte} 
Let $(M,H,g^M)$ be a smooth sub-Riemannian contact manifold with $C^k$- bounded geometry ($k\in\N$),  and 
let $\{\chi_j\}$ be a partition of the unity subordinated to the atlas $\mc U:=\{ B(x_j, \rho), \phi_{x_j}\}$ of Lemma \ref{covering}. 
 We stress explicitly that $\phi_{x_j}^{-1}(\mathrm{supp}\;\chi_j) \subset B(e,1)$.
If $\alpha$ is a Rumin's differential form on $M$, we say that $\alpha\in W^{\ell,p}_{\mc U}(M, E_0^\bullet)$ 
for $\ell=0,1,\cdots, k-1$,  and $p\ge 1$ if
$$
\phi_{x_j}^\sharp (\chi_j\alpha) \in W^{\ell,p}(\he n, E_0^\bullet)\quad\mbox{for $j\in \mathbb N$}
$$
(notice that 
$\phi_{x_j}^\sharp (\chi_j\alpha)$ 
is compactly supported in $B(e,1)$ and therefore can be continued by zero on
all $\he n$). Then we set
\begin{equation}
	\label{spazi-Sobolev-M}
\|\alpha\|_{W^{\ell,p}_{\mc U}(M, E_0^\bullet)}:= \left(\sum_j \|\phi_{x_j}^\sharp (\chi_j\alpha)\|^p_{W^{\ell,p}(\he n, E_0^\bullet)}\right)^{1/p}.
\end{equation}
\end{definition}
An other uniform covering and other choices of controlled charts lead to an equivalent norm,
so the definition of the Sobolev spaces ${W^{\ell,p}_{\mc U}(M, E_0^\bullet)}$ 
is indeed 
not depending on the atlas $\mc U$, as shown in 
\cite{BFP2}, Proposition 4.13. 
Therefore, from now on  we drop the index $\mc U$ from the notation of Sobolev norms and we shall write simply $
W^{\ell,p}(M, E_0^\bullet).$

{\begin{remark}\label{norme_lp}
Notice that if $\ell=0$ in definition \eqref{spazi-Sobolev-M}, the norm in $W^{0,p}(M, E_0^\bullet)$ is equivalent to the norm $L^p(M, E_0^\bullet)$ associated to the volume form $\mu:=\theta^M\wedge (d \theta^M)^n$ (see the Introduction).  

\end{remark}

\begin{proof}
Let us denote the norm in $W^{0,p}(M, E_0^\bullet)$ by
 \begin{equation*}
	\label{spazi-lp-M}
|||\alpha|||_{L^p(M, E_0^\bullet)}:=\|\alpha\|_{W^{0,p}(M, E_0^\bullet)}= \left(\sum_j \|\phi_{x_j}^\sharp (\chi_j\alpha)\|^p_{L^{p}(\he n, E_0^\bullet)}\right)^{1/p},
\end{equation*}
and by 
$$\|\alpha\|_{L^p_\mu(M,E_0^\bullet)}$$
the $L^p$ norm associated with $\mu$.

First of all, since $\chi_j\alpha$ is compactly supported and $\phi_{x_j}^\sharp(\theta^M\wedge (d\theta^M)^n)=\phi_{x_j}^\sharp(\theta^M)\wedge \phi_{x_j}^\sharp((d\theta^M)^n) =\theta^{\mathbb H}\wedge (d\theta^{\mathbb H})^n$, we have  $\|\phi_{x_j}^\sharp (\chi_j\alpha)\|^p_{L^{p}(\he n, E_0^\bullet)}\approx  \|\chi_j\alpha\|^p_{L^{p}_{\mu }(M, E_0^\bullet)}$ hence
$$
|||\alpha|||_{L^p(M, E_0^\bullet)}=\left(\sum_j \|\phi_{x_j}^\sharp (\chi_j\alpha)\|^p_{L^{p}(\he n, E_0^\bullet)}\right)^{1/p}\approx \left(\sum_j \|\chi_j\alpha\|^p_{L^{p}_{\mu }(M, E_0^\bullet)}\right)^{1/p}\,.
$$

Hence we are left to show that 
$$
\left(\sum_j \|\chi_j\alpha\|^p_{L^{p}_{\mu }(M, E_0^\bullet)}\right)^{1/p}\approx \|\alpha\|_{L^{p}_{\mu }(M, E_0^\bullet)}\,.
$$

To this aim, first notice that, since $\{\chi_j\} $ is a partition of the unity $|\chi_j|\le 1$, 
hence $\sum_j |\chi_j|^p\le \sum_j |\chi_j|\le 1$ and therefore 

$$
\sum_j\int_M |\chi_j\alpha|^p d\mu=\int_M \sum_j|\chi_j\alpha|^p d\mu =\int_M \sum_j|\chi_j|^p|\alpha|^p\le \|\alpha\|^p_{L^p_\mu(M,E_0^\bullet)}.
$$

On the other hand,  since $\alpha=\sum_j \chi_j\alpha$,
$$
\|\alpha\|^p_{L^p_\mu(M,E_0^\bullet)}\le  \|\sum_j\chi_j\alpha\|^p_{L^p_\mu(M,E_0^\bullet)}= \int_M|\sum_j\chi_j|^p|\alpha|^pd\mu.
$$
Now, for any $x\in  M$, $\sum_j\chi_j(x)$ is a finite sum with a number of terms less than or equal to $N$ (see Lemma \ref{covering}-iii)) hence there exists a constant $c_N$ so that $|\sum_j\chi_j|^p\le c_N \sum_j|\chi_j|^p$ and in the inequality above we get

$$
\|\alpha\|^p_{L^p_\mu(M,E_0^\bullet)}\le c_N\sum_j\int_M|\chi_j|^p|\alpha|^pd\mu=c_N\sum_j \|\chi_j\alpha\|^p_{L^{p}_{\mu }(M, E_0^\bullet)}.
$$

\end{proof}}

\medskip

Thanks to the previous remark,  from now on we shall denote with the same symbol $\|\,\cdot\,\|_{L^p(M, E_0^\bullet)}$ the two equivalent norms $
\|\,\cdot\,\|_{L^p_\mu(M,E_0^\bullet)}$ and $
|||\,\cdot\,|||_{L^p(M, E_0^\bullet)}$.

In the sequel of the paper we will need a covering of $M$ with balls of suitable, fixed radius, which has the same properties as the covering of the Lemma \ref{covering}. In analogy with what happens in the Riemannian setting (see e.g. \cite{triebel}, 7.2.1) we have:
\begin{remark}\label{covering ad hoc}
If $\eta >0$ is small, then there exists   an at most countable covering  of $M$ with balls  $\{ B(a_\ell, \eta)\}$ which satisfy the same properties as the covering given in Lemma \ref{covering}, in particular the covering is uniformly locally finite, i.e. there exists $N=N(\eta)\in \mathbb N$ such that any point
of $M$ has an open neighborhood of radius $\eta$ that is covered by at most $N(\eta)$ balls of the covering. Moreover the functions of the atlas $\left\{B(a_\ell,\eta), \phi_\ell\right\}$, with $\phi_\ell: B(e, \eta/C_M)\to M$  satisfy conditions (2) and (3) of Definition \ref{contact bis}  with  constants depending only on $C_M$ but  independent of $a_\ell$.
In particular,  we notice that, if $\alpha$ is supported in $\phi_\ell(B(e,\eta/C_M))$, then by the Definition
\ref{contact bis} the norms
$$
\| \alpha\|_{W^{m,p}(M,E_0^\bullet)} \qquad \mbox{and}\qquad \|\phi_\ell^\sharp \alpha\|_{W^{m,p}(\he n,E_0^\bullet)}
$$
are equivalent, with equivalence constants independent of $\ell$.
\end{remark}

\begin{proof}

%

For any $x\in M$, let 	$\phi_x\colon B(e,1)\subset \mathbb{H}^n \to M$ be a map satisfying the conditions contained in Theorem \ref{delirio}.

Let $\{B({x_j},\rho)\}$ be a countable locally finite subcovering of $\{\phi_{x}(B(e,1))\, , x\in M\}$ as in Lemma \ref{covering}. 
and let $\phi_{x_j}$ be the corresponding bounded contact charts from the unit Heisenberg ball i.e. $\phi_{x_j}: B(e,1)\to B({x_j},\rho)$. We show now that any ball $B({x_j},\rho)$
 has a finite subcover of balls of radii $\eta$.

 Without loss of generality, we can assume that the $\phi_{x_j}$ are defined on a larger ball $B(e,\lambda)$, where $\lambda >1$ is fixed and they still satisfy conditions (2) and (3) of Definition \ref{contact bis}.

Let $\eta>0$. We can cover the ball $B(e,1)$ by a finite number $k=k(\eta)$ of balls of radii $\frac{\eta}{C_M}$   
\begin{equation}\label{palline}
	B(e,1)\subset\bigcup_{i=1}^{k}B(z_i, \frac{\eta}{C_M})\subset B(e,\lambda)\,.
\end{equation}

For any $z_i$ as in \eqref{palline},  we define the map $\phi_{x_j}^i: B(e,\frac{\eta}{C_M})\to M$ as follows
$$
\phi_{x_j}^i(z):=\phi_{x_j}\circ \tau_{z_i}(z)=\phi_{x_j}(z_i\cdot z)
$$
Notice that $
\phi_{x_j}^i(e)=\phi_{x_j}(z_i\cdot e)=:a_{x_j}^i\in \phi_{{x_j}}(B(e,1))$ 

Since the map $\phi_{x_j}$ is a contact map, and we have  composed with a translation, the map $
\phi^i_{x_j} $ is a smooth contact map. Indeed,  since $\phi_{x_j}$ is a contact map we have $({\phi_{x_j}})^\sharp\theta^M=\theta^{\he {}}$, hence

$$
({\phi_{x_j}^i})^{\sharp }\theta^M=(\phi_{x_j}\circ \tau_{z_i})^\sharp\theta^M= \tau_{z_i}^\sharp\circ ({\phi_{x_j}})^\sharp\theta^M=\tau_{z_i}^\sharp\theta^{\he {}}=\theta^{\he {}}\,.
$$

Moreover the maps  $\phi_{x_j}^i$ satisfies conditions (2) and (3) of Definition \ref{contact bis} with  constants depending only on $C_M$ but  independent of $x_j$. We can see for example that (2) still holds for  $\phi_{x_j}^i$, as a consequence of the left-invariance of the Kor\'any distance \eqref{def_distance}. Indeed, for all $p,q\in B(e,1)$ we have
 \begin{equation}\label{sopra}
d_M(\phi_{x_j}^i (p), \phi_{x_j}^i(q)) =d_M(\phi_{x_j} (z_i\cdot p), \phi_{x_j}(z_i\cdot q))\approx d(z_i\cdot p,z_i\cdot q)=d(p,q)\, , 
\end{equation}
where the symbol $\approx$ means that we have used the same constant $C_M$  given in \eqref{bilip}.

Now, reasoning as in \eqref{sopra}, for any $i=1,\cdots, k=k({\eta})$ the balls $B(a_{x_j}^i, \eta )\supset \phi_{x_j}^i(B(e,\frac{\eta}{C_M}))$  hence $$\bigcup_{i=1}^{k({\eta})} B(a_{x_j}^i, \eta )$$ is a finite cover of the set $B({x_j},\rho)$.


Therefore  any ball $B({x_j},\rho)$
of the countable covering $\left\{B(x_j,\rho)\right\}$ of $M$ has a finite subcover of balls of radii $\eta$ and eventually
$\left\{B(a_{x_j}^i,\eta)\right\}$ is an at most countable covering of $M$ uniformly locally finite. 

From now on we shall denote the covering  $\left\{B(a_{x_j}^i,\eta), \phi_{x_j}^i\right\}$ simply
by $\left\{B(a_\ell,\eta), \phi_\ell\right\}$. By construction any point
of $M$ has an open neighborhood of radius $\eta$ that is covered by at most $N(\eta)$ balls of the covering.
We notice that the number $N(\eta)\approx k(\eta)N(M)$ where the symbol $\approx$ means that there are constants depending only on the geometry of $M$ (i.e. on $C_M$ and $r$) and $N(M) $ is the number appearing Lemma \ref{covering}. 


\end{proof}

\begin{remark}
\label{eta tilde} We notice that, if we use a covering of $M$ with balls of radius $\eta$ small as above,  the constant $c_N$ which gives the equivalence between the two norms in Remark \ref{norme_lp}, will depend also on $\eta$.
\end{remark}

\subsection{Symplectic basis and orthogonal linear transformations}

In the sequel of the paper we need to cover $M$ with atlases that enjoy further properties besides those contained in Definition \ref{contact bis}. First of all, in the next theorem we observe that we can replace the contactomorphism  $\phi_x$ appearing in Definition \ref{contact bis} with another  contactomorphism which ``sends'' an orthonormal symplectic basis of $\ker\theta^M_x$ into the  canonical orthonormal symplectic basis of $\ker \theta^{\mathbb H}_e$ (see \eqref{campi W}) and still depending only on the bounded geometric constants and not on the point $x$.  We begin with the following remark.

\begin{remark}\label{mah}
Given a contact manifold $M$,  there exists an orthonormal  basis of $\ker\theta^M$ pointwise.
\end{remark}
\begin{proof}

 This can be shown by simply considering the endomorphism $J: \ker\theta^M\to \ker\theta^M$, with $J^2=-Id$. If we follow the steps of the proof of Theorem 2.1.3 in \cite{mcduff_salamon}, and choose $Z_1\in \ker\theta^M$ to be a unit vector field, i.e. $$1=g^M(Z_1,Z_1)=d\theta^M(Z_1, JZ_1),$$ then also $JZ_1$ is a unit vector field, as 
\begin{align*}
g^M(JZ_1,JZ_1)&=d\theta^M(JZ_1, J^2Z_1)=-d\theta^M(J^2Z_1, JZ_1)\\
&=-d\theta^M(-Z_1, JZ_1)=g^M(Z_1,Z_1)=1\,.
\end{align*}
 
Notice also that $g^M(Z_1, JZ_1)=d\theta^M(Z_1, J^2Z_1)=-d\theta^M(Z_1, Z_1)=0$, so $Z_1$ and $JZ_1$ are orthonormal.

In order to extend $\lbrace Z_1,JZ_1\rbrace$ to an orthonormal basis of $\ker\theta^M$, let us consider a
unitary vector field $Z_2\in \mathrm{span} \{ Z_1, J Z_1\}^\perp \cap \ker \theta^M$ . Arguing as above, $g^M(Z_2,J Z_2)=0$. To prove that the vectors fields $Z_1, JZ_1, Z_2,JZ_2$ are all orthonormal, we are left to show that  $g^M(Z_1,JZ_2)=0 $ and $g^M(JZ_1,JZ_2)=0$. Indeed,
$$g^M(Z_1, JZ_2)=d\theta^M(Z_1, J^2Z_2)=-d\theta^M(JZ_1, JZ_2)=-g^M(JZ_1,Z_2)=0$$ by construction, and analogously $$g^M(JZ_1,JZ_2)=d\theta^M(JZ_1,J^2Z_2)=d\theta^M(Z_1,JZ_2)=g^M(Z_1,Z_2)=0$$ by construction.

Likewise, one can repeat the same reasoning $n$ times and construct an orthonormal symplectic basis for $\ker\theta^M$.

\end{proof}

From now on, the basis $\{Z_1,\cdots,Z_n,JZ_1,\cdots, JZ_n\}$ will be denoted by $$\{W_1^M,\dots, W_{2n}^M\}$$ and we will refer to it as an  orthonormal symplectic basis of $\ker\theta^M$.

We   recall now the following definition. 


	\begin{definition}\label{hodge}
		Let $V$ and $W$ be real vector spaces of dimension $N$, both endowed with scalar products $\langle\cdot,\cdot\rangle_V$ and $\langle\cdot,\cdot\rangle_W$ respectively. We say that the linear map
		\begin{align*}
			T\colon V\to W
		\end{align*}
		is an {orthogonal linear transformation} if, given $\lbrace e_1,\ldots,e_N\rbrace$ an orthonormal basis of $V$, then $\lbrace Te_1,\ldots,Te_N\rbrace=\lbrace \varepsilon_1,\ldots,\varepsilon_N\rbrace$ is an orthonormal basis of $W$.
	\end{definition}

The following proposition
   follows easily from a result of \cite{agrachev_barilari_boscain}.

\bigskip
	\begin{proposition}\label{bf} Let $a\in M$ be a fixed point, and $W_1^M,\dots, W_{2n}^M$  an orthonormal symplectic
basis of   $\mathrm{ker}\,\theta^M$  in a neighboorood of $a$. Then there exist $\epsilon=\epsilon(a)>0$ and  a smooth family of horizontal curves  $\gamma^\ell: [0,\epsilon]\rightarrow M$,  for $\ell=1,\cdots,2n$, such that
 \begin{itemize}
 \item[i)] $\gamma^\ell (0)=a$;
 \item [ii)] $(\gamma^\ell)'(0) = W_{\ell}^M(a)$ ;
 \item [iii)] $
 d_M(\gamma^\ell(t),a)=\int_0^t g^M((\gamma^\ell)'(s), (\gamma^\ell)'(s))_{\gamma^\ell(s)}^{1/2} ds.
 $
 
 \end{itemize}

\end{proposition}

\begin{proof} Fix $\ell\in\{1,\cdots,2n\}$. Following \cite{agrachev_barilari_boscain}, Section 4.3.1, denote by $H:T^*M\to \mathbb R$
the sub-Riemannian Hamiltonian associated with $(M,\ker\theta^M,g^M)$, and let $\lambda: [0,1]\to T^*M$
the normal extremal, i.e. the solution of
$$
\lambda'(t) = \vec{H}(\lambda(t)),
$$
with $\lambda(0)= (W^M_{\ell}(a))^\natural$ (recall, the Hamiltonian vector field $\vec{H}$ is defined by $d\theta^M(\vec{H},\, \cdot\,)=-dH$) , so that
$$
\Scal{\lambda(0)}{W^M_{j}(a)}=\delta_{\ell j}.
$$

Then the assertion follows from Remark 4.28 and Theorem 4.65 in \cite{agrachev_barilari_boscain}.

\end{proof}

The main result of this section is the following theorem.

\begin{theorem}\label{delirio}
 There exist $0<r'<r$,  and $0<\mu\le 1$ (all depending only on the bounded geometry constants) such that for any $a\in M$, if we denote by $\left\{W_{1,a}^M,\ldots,W_{2n,a}^M\right\}$  the orthonormal symplectic  basis of $\ker\theta^M_a$ (and, as in   \eqref{campi W}, $\{ W_{1,e}^\mathbb{H},\ldots,W_{2n,e}^\mathbb{H}\}$ is the orthonormal symplectic basis of $\ker\theta^\mathbb H_e$),  there exists a contact map $\psi_a$ (that is $\psi_a^\sharp(\theta^M)=\theta^\mathbb{H}$),
\begin{align*}
			\psi_a\colon B(e,\mu)\subset \mathbb{H}^n \to M
		\end{align*}
 satisfiying $B(a,r')\subset \psi_a(B(e,\mu))$ and conditions (2) and (3) given in Definition \ref{contact bis} and  
 such that
		\begin{itemize}
			\item[i.] $\psi_a(e)=a$ 
			\item[ii.] $(d\psi_a)_e W_{j,e}^\mathbb{H}=W_{j,a}^{M}$ for $j=1,\ldots,2n$,
			and $(d\psi_a)_e\xi_e^\mathbb{H}=\xi_{a}^M$. 
			In particular, 
			the map 
			\begin{align*}
				(d\psi_a)_e\colon T_e\mathbb{H}^n\to T_aM\\
				v\mapsto (d\psi_a)_e(v)
			\end{align*}
			is an orthogonal linear map.
		\end{itemize}
	\end{theorem}
	
	\begin{proof}
		 {Let $a\in M$ and let}  $\phi_a : B(e,1)\to M$ be  a  contactomorphism
		as in Definition \ref{contact bis} i.e.
\begin{enumerate}
  \item $B(a,r)\subset\phi_a(B(e,1))$;
  \item
  \begin{equation}\label{bilip_due}
\frac{1}{C_M}d(p,q)\le  d_M(\phi_a (p), \phi_a(q)) \le C_Md(p,q)\qquad \mbox{for all $p,q\in B(e,1)$};
\end{equation}
  \item coordinate changes $\phi_b^{-1}\circ \phi_a$ and their first $k$ derivatives with respect to unit left-invariant horizontal vector fields are bounded by $C_M$.
\end{enumerate}
\indent We can also assume that $\phi_a (e)=a$.
	
		We consider the map
		\begin{align*}
			(d\phi_a)_e\colon T_e\mathbb{H}^n\to T_aM\,.
	\end{align*}

The map $\phi_a^{-1} : \phi_a ( B(e,1))\to  B(e,1)$ defines, by pushforward,
the vector fields
$$
\left\{\hat W_1,\cdots , \hat W_{2n}\right\}:=\big\{d\phi_a^{-1}(W_1^M),\cdots d\phi_a^{-1}(W_{2n}^M)\big\}
$$
that are a symplectic basis of  $\ker\theta^\mathbb H$. 
Indeed, since $\phi_a$ is a contact map i.e. $$\phi_a^\sharp(\theta^M)=\theta^\mathbb{H},$$ for example if $j=i+n$  we have

\begin{equation*}
\begin{split}
d\theta^\mathbb H&\Big(\hat W_i, \hat W_{j}\Big)=d\theta^\mathbb H\Big(d\phi_a^{-1}(W_i^M),\ d\phi_a^{-1})(W_j^M)\Big)
\\&
=d(\phi_a^\sharp\theta^M)\Big(d\phi_a^{-1}(W_i^M),\ d\phi_a^{-1}(W_j^M)\Big)
\\&=\phi_a^\sharp(d\theta^M)\Big(d\phi_a^{-1}(W_i^M),\ d\phi_a^{-1}(W_j^M)\Big)
\\&=\Scal{\ d\theta^M\ }{\  d\phi_a\,d\phi_a^{-1}(W_i^M),d\phi_a\,d\phi_a^{-1}(W_j^M)\ }=d\theta^M\Big( W_i^M,W_j^M\Big)=\delta_{ij}
\end{split}
\end{equation*}

Then in particular $(\hat W_1(e),\cdots , \hat W_{2n}(e))$
 can be identified with a symplectic basis of $\mathbb{R}^{2n}$.
Hence, if we denote by $
(e_1,\cdots , e_{2n})$ the canonical basis of $\mathbb{R}^{2n}$,  there exists a matrix $\tilde A=\tilde A_{\phi_a}\in Sp(2n)$ such that 
$$
\tilde A W_{i,e}^\mathbb{H}=\tilde A e_i=\hat W_i(e)\quad \quad i=1,\cdots,2n\,,
$$
%
%
(we stress that the matrix $\tilde A$ depends on  $d\phi_a^{-1}$).
	It is well defined the (Euclidean) linear contact map $L:\he n\to \he n$ associated with the matrix $$A:=\left( 
	\begin{aligned}
	\tilde A \quad \quad& 0_{2n\times 1}\\
	0_{1\times 2n } \quad \quad	& 1
	\end{aligned}
	\right)$$
	with $A\in GL(\R^{2n+1},\R^{2n+1})$. In particular, $A$ induces an automorphism of the group.
	
	\bigskip
	
	We are now ready to show that there exist $0<r'<r$ and $0<\mu\le 1$ (depending only on the constant $C_M$ appearing in \eqref{bilip_due}) and a contact map $\psi_a$ that satisfies the property of Definition \ref{contact bis}, with constants depending only on the bounded geometry, but independent of $a$.
	
	{\it Claim:} We claim that the norm of $ \hat W_i(e)$  can be bounded from above and below by constants depending only on the constant $C_M$ appearing in Definition \ref{contact bis} and not on the point $a$. We can then write  $\|\hat W_i(e)\|\approx 1$ (independently of the point $a\in M$), where the symbol $\approx$ means that the constants appearing above depend only on $C_M$.
	
	\medskip
	
	Let assume for a while that the claim is true. It follows that also the 
	norm of the matrix $\tilde A$ is controlled from below and above by a constants depending only on the constant $C_M$ appearing in Definition \ref{contact bis}, i.e. $1/C_M\le\|\tilde A\|\le C_M$. 
	
	\medskip


	%

First of all, we notice that if $p=(p_1,\cdots, p_{2n}, p_{2n+1})=:(p',p_{2n+1})\in \he n$ and $L(p)=A(p', p_{2n+1})=(\tilde A p', p_{2n+1})$, then ($C_M>1$)
	\begin{equation}\label{autovalori}
\begin{split}	
C_M^{-4}d(p,e)^4&\le  C_M^{-4}(|\tilde A^{-1}\tilde A p'|^4+16 p_{2n+1}^2)
\\
&\le (|\tilde Ap'|^4+16p_{2n+1}^2)=	d(L(p),e)^4\le  C_M^4d(p,e)^4\,.
\end{split}
	\end{equation} 
	Therefore,
	if $0<\mu\le 1$,  we have 
\begin{equation}\label{palle}
\begin{split}
	 &L(B(e,\mu))\subset B(e,1) \,
	\end{split}
\end{equation}
if we take $\mu\le 1/C_M$.
	

	Hence, we can define
	$$
	\psi_a:=\phi_a\circ L_{|B(e,\mu)}:  B(e,\mu)\to M\,.
	$$
	
	The map $\psi_a $ is a contactomorphism   and by construction satisfies\begin{align*}
			W_{j,a}^M=(d\psi_a)_{e}W_{j,e}^\mathbb H\,.
		\end{align*}
	
	We show now that if $r'<\frac{\mu}{ C_M}$  we have 
	\begin{equation}\label{erre}
		B(a,r')\subset\psi_a(B(e,\mu))
	\end{equation}
	Indeed, let $b\in B(a,r')$. There exists $p\in B(e,1)$ such that $b=\phi_a(p)$ (if $r'<r$). By \eqref{bilip_due}, if we take $r'<\frac{\mu}{C_M^2}$ we have
	$$
	d(p,e)\le C_M d_M(\phi_a(p),a)=C_M d_M(b,a)<{\frac{\mu}{C_M}}
	$$
	i.e.  $p\in B(e,\frac{{\mu}}{C_M})$. Since $L$ is a linear isomorphism there exists $q\in\he n$ such that $p=L(q)$. We show that $q\in B(e,\mu)$. Indeed, by \eqref{autovalori}, $d(q,e)\le C_Md(L(q), e)=C_Md(p,e)<C_M \frac{{\mu}}{C_M}=\mu$. Hence, 
	$b=\phi_a(p)=(\phi_a\circ L)(q)=\psi_a(q)\in \psi_a(B(e,\mu))$ and  \eqref{erre} holds.

	We need to show now that if $p,q \in B(e,\mu)$ then the condition \eqref{bilip} in Definition \ref{contact bis} is satisfied i.e.
	$$
	d_M(\psi_a(p),\psi_a(q))\approx d(p,q)
	$$
	where, here, the symbol $\approx$ means that the constants appearing above depend only on $C_M$. Notice that if $p,q \in B(e,\mu)$ then $L(p),L(q) \in B(e,1)$, by \eqref{palle}, and
	\begin{equation*}
	\begin{split}
	d_M(\psi_a(p),\psi_a(q))&=d_M(\phi_a(L(p)), \phi_a(L(q)))\overbrace{\approx}^{by \, \eqref{bilip_due}  } d(L(p),L(q))
	\\
	&=\rho(L(p)^{-1}\cdot \ L(q))=\rho(L(p^{-1}\cdot\ q))=d(L(p^{-1}\cdot\ q),e)\,.
	\end{split}
	\end{equation*}
	
As in \eqref{autovalori},
$$d(L(p^{-1}\cdot\ q),e)\approx \ d(p, q),$$
		hence 
		$$  d_M(\psi_a(p),\psi_a(q))\approx  \ d(p,q)$$
		where, again,  the equivalence constants depend only on the constant $C_M$ appearing in Definition \ref{contact bis} and hence independent of $a\in M$.
		
		%
		
	{	We have now to check that the map $\psi_a$ satisfies also  the third requirement of Definition \ref{contact bis}. Let $a,b\in M$, and consider the maps 
	 $\psi_a=\phi_a\circ L$ and $\psi_{b}=\phi_b\circ \hat L$ constructed as above. We remind that the maps $L$ and $\hat L$, being linear contact maps,  preserve horizontal derivatives. Since, by (3) in Definition \ref{contact bis}, coordinate changes $\phi_b^{-1}\circ \phi_a$ and all their first $k$ derivatives with respect to unit left-invariant horizontal vector fields are bounded by $C_M$, it follows that $\psi_b^{-1}\circ\psi_a=\hat L^{-1}\circ \phi^{-1}_b\circ \phi_a\circ L$ enjoy the same property.
	}
		%
		
		\bigskip
		
		{\it We are left with the proof of the Claim}. 
		
	Let $W_1^M,\dots, W_{2n}^M$ be an orthonormal symplectic
basis of   $\mathrm{ker}\,\theta^M$  in a neighboorood of $a$. 	By Proposition \ref{bf}, for any $j=1,\cdots,2n$ there exist  a curve $\gamma^j(t)$  for $t$ small, such that:
 $$d_M(\gamma^j(t),a)=\int_0^t g^M((\gamma^j)'(s), (\gamma^j)'(s))_{\gamma^j(s)}^{1/2}ds .$$
Notice that the basis is orthonormal, hence for any $j$, $g^M((\gamma^j)'(0), (\gamma^j)'(0))_{a}^{1/2}=g^M(W_{j,a}^M, W_{j,a}^M)_{a}^{1/2}=1$. 
 
Let us take  the map $\phi_a$ 	considered at the very beginning of the proof. Hence, again by condition	\eqref{bilip_due},   for	$t\neq 0$ small, we have
		\begin{equation}\label{giovedi}
 \frac{1}{t}\int_0^t g^M((\gamma^j)'(s), (\gamma^j)'(s))_{\gamma^j(s)}^{1/2}ds=\frac{d_M(\gamma^j(t),a)}{t}\approx   \frac{d(\phi_a^{-1} \gamma^j(t),e)}{t},
 \end{equation}
		where the symbol $\approx$ means that the constants appearing  depend only on $C_M$.
		
Set
$
\sigma^j(t):=\phi_a^{-1}\gamma^j(t)$, then $\sigma^j(0)=\phi_a^{-1} \gamma^j(0)=e$.  Since  $\hat W_j=d\phi_a^{-1}(W_j^M),
$
we have 
 $(\sigma^j)'(0)=(d\phi_a^{-1})_e(\gamma^j)'(0)=(d\phi_a^{-1})_eW_{j,a}^M=\hat W_j(e)$.

The map $\phi_a$ is a contact map, then the  vector fields $\hat W_j$ are horizontal vector fields.

From now on, we argue with a fixed vector field $\hat W_j$ hence,  for sake of simplicity, we shall drop the index $j$ writing $\hat W$ instead of $\hat W_j$ and $\sigma(t)$ instead of $\sigma^j(t)$. Hence we have
\begin{equation}\label{campi sigma}
\hat W(p)=\sum_{k=1}^{2n}\lambda_k(p)W_{k,p}^\mathbb{H}
\end{equation}
and 
$\sigma(t)=(\sigma_1(t),\cdots, \sigma_{2n}(t),\sigma_{2n+1}(t))$ satisfies:

$\sigma_k'(t)=\lambda_k(\sigma(t))$ if $1\le k\le 2n$, 


$\sigma_{2n+1}'(t)
= \frac12\sum_{k=1}^n\left(\sigma_k'(t) \sigma_{k+n}(t)-\sigma_{k+n}'(t) \sigma_{k}(t)\right)$.

\medskip
		
		By Taylor's formula we have, for $t\to 0$:
		
		$\sigma_k(t)=t\lambda_k(\sigma(0))+\ O(t^2)=t\lambda_k(e)+\ O(t^2)$ if $1\le k\le 2n$, and hence also $\sigma'_k(t)=\lambda_k(e)+\ O(t)$ if $1\le k\le 2n$.		
		Replacing these expressions in $\sigma_{2n+1}'(t)$ we get
	$$	
	\begin{aligned}
		\sigma_{2n+1}'(t)
&= \frac12\sum_{k=1}^n\left((\lambda_k(e)+O(t))(t\lambda_{k+n}(e)+O(t^2))-(\lambda_{k+n}(e)+O(t))(t\lambda_{k}(e)+O(t^2)\right)
\\
&=O(t^2)=o(t^3)\,.
\end{aligned}	
$$	
Therefore, for $t\to 0$, 
		$$
		\begin{aligned}
		\sigma(t)&=\left(t(\lambda_1(e)+\ o(1)),\cdots, t(\lambda_{2n}(e)+\ o(1)), t^2\cdot o(1)\right)\\
		&=\delta_t\left(\lambda_1(e)+\ o(1),\cdots, \lambda_{2n}(e)+\ o(1),  o(1)\right)\,,
		\end{aligned}
		$$
		where $\delta_t$ is the dilation defined in \eqref{dilatazione}.
		If we take $p=e$ in \eqref{campi sigma}, $\hat W(e)=\left(\lambda_1(e),\cdots, \lambda_{2n}(e), 0\right)$, hence 
		$$
		\sigma(t)=\delta_t\left(\hat W(e)+o(1) \right)\,
		$$
		Thus $$\frac{d(\sigma(t),e)}{t}=\frac{\rho(\sigma(t))}{t}= \|\hat W_j(e)+o(1)\|\,.$$
		Therefore  by \eqref{giovedi} it holds
		$$
		\frac{1}{C_M} \|\hat W_j(e)+o(1)\|\le \frac{1}{t}\int_0^t g^M((\gamma^j)'(s), (\gamma^j)'(s))_{\gamma^j(s)}^{1/2}\le {C_M} \|\hat W_j(e)+o(1)\|\,,
		$$
		as $t\to 0$. 
		Passing now to the the limit for $t\to 0$ , and keeping in mind that $g^M(\gamma'(0), \gamma'(0))_{a}^{1/2}=1$ we get

		$$
		\frac{1}{C_M} \|\hat W_j(e)\|\le 1 \le C_M \|\hat W_j(e)\|,$$ and the claim is proved.
		
	\end{proof}
	
	\begin{remark}
	The new constants $r'$ and $\mu$ depend on $r, \, C_M, \ 1$.  To avoid cumbersome notation we shall  still denote these constants, { respectively}, by $r$, $C_M$  and $1$.  When we want to cover $M$ with balls of radius $\eta$ (small), in analogy with the previous proof, the constants $r'$ and $\mu$ will depend also on $\eta$. From now on, for sake of simplicity, we will denote with the same constant $C_M$ any dependence on $r, \, C_M$.
	\end{remark}
	


%
%
%
%
%

%
%
%

\section{Some results about the lack of commutation between the pull-back and the co-differential}\label{conticodiff}

In this section we take the map $\psi: B(e,1)\to M$ as in Theorem \ref{delirio}. 

Let $\alpha$ be a smooth differential form on $M$ and set $\beta:=\psi^\sharp\alpha$ the pullback of $\alpha$. As already pointed out in Proposition \ref{memory},  pulling back $d_c^M\alpha$ gives $d_c^{\mathbb H}\beta$. On the contrary, since $\delta_c^M\alpha$ involves the Hodge-$\ast$ operator, it turn out that its pullback is not $\delta_c^{\mathbb H}\beta$, i.e.
 $\psi^\sharp\delta_c^M \neq \delta_c^{\mathbb H}\psi^\sharp $. In this section we shall examine the  relation between $\psi^\sharp\delta_c^M $ and $\delta_c^{\mathbb H}\psi^\sharp $. Remeber that both $\delta_c^{\mathbb H}$ and $\delta_c^{M}$  are equal to $\pm\ast d_c^{\mathbb H}\ast$ and $\pm\ast d_c^{M}\ast$ respectively. In the sequel we shall always drop the sign since   we are only interested in estimates of norms.

In Definition \ref{pullback} and  Proposition \ref{giaquinta} below we recall  some preliminary  notations and results(see  e.g. \cite{GMS}, Section 2.1).

{
	
		Let $V$ and $W$ be real vector spaces of dimension $N$, both endowed with scalar products $\langle\cdot,\cdot\rangle_V$ and $\langle\cdot,\cdot\rangle_W$ respectively.
	
	
	If we denote by $\ast$ the Hodge-$\ast$ operator in $V$, the following equality holds
	\begin{align*}
		u\wedge \ast v:=\langle u,v\rangle_V \ e_1\wedge\cdots\wedge e_N.
	\end{align*}
	where $\{e_1,\cdots,e_N\}$ is an orthogonal basis of $V$.
	The Hodge-$\ast$ operator on $W$ is defined analogously. 
	
	To fix our notations, we remind the following definition (see e.g. \cite{GMS}, Section 2.1).
\begin{definition}\label{pullback} If $V,W$ are finite dimensional linear vector spaces and
$L:V\to W$ is a linear map, we define
$$\Lambda_{h} L :
\bigwedge\nolimits_{h} V\to \bigwedge\nolimits_{h} 
W$$ 
as the linear map defined by
$$
(\Lambda_{h} L) (v_1\wedge\cdots\wedge v_h)=
L(v_1)\wedge\cdots\wedge L(v_h)
$$
for any simple $h$-vector
$v_1\wedge\cdots\wedge v_h\in
\bigwedge\nolimits_{h} V$,  and 
$$\Lambda^{h} L :
\bigwedge\nolimits^{h} W\to \bigwedge\nolimits^{h} V$$
as the linear map defined by
$$
\Scal {(\Lambda^h L)(\alpha)}{v_1\wedge\cdots\wedge v_h}=
\Scal{\alpha}{(\Lambda_hL)(v_1\wedge\cdots\wedge v_h)}
$$
for any $\alpha\in \bigwedge\nolimits^{h} W$
 and any simple $h$-vector
$v_1\wedge\cdots\wedge v_h\in
\bigwedge\nolimits_{h} V$.
\end{definition}

	\begin{proposition}\label{giaquinta} If $V$ and $W$ are $N$-dimensional vector spaces and $L:V\to W$
		is an orthogonal linear transformation, we have:
		$$
		\ast (\Lambda_h L )=  (\Lambda_{N-h}L) \ast.
		$$
	\end{proposition}
	
	\bigskip
	
	%
%
	
 %
%
%
%
%
%

	\bigskip
	
	Let $\{W_1^M, \dots, W_{2n}^M\}$ 
	be  an orthonormal symplectic
	basis of $\ker\theta^M$, if we
	 denote by $\xi^M$ the Reeb vector field on $M$ then the metric  can be
	extended to a Riemannian metric on $TM$ (still denoted by $g^M$), so that
	$\{W_1^M, \dots, W_{2n}^M,\xi^M\}$ is an orthonormal frame of $TM$.
	}
	Take $W_1^{\he{}}, \dots, W_{2n}^{\he{}}$
	 the standard orthonormal symplectic basis of $\ker\theta^{\he{}}$ (see \eqref{campi W}).

	If $m\in M$, we denote by $$k_m : T_mM \to \mathbb R^{2n+1}$$ the map which associates
	to a vector ${v}\in T_mM$ its coordinates with respect to the basis
	$\{W_{1,m}^M, \dots, W_{2n,m}^M,\xi^M_m\,\}$, and by
	$$f_x : T_{x}\he n  \to \mathbb R^{2n+1}$$ the analogous map which associates
	to a vector ${v}\in T_{x}\he n $ its coordinates with respect to the basis
	$\{W_{1,x}^{\he{}}, \dots, W_{2n,x}^{\he{}},W^{\he{}}_{2n+1,x}\}$.
	
	We have the following property.
	
	\begin{lemma}\label{RMK MAPPA}
		If $\psi: B(e,1)\to M$ is a map as in Theorem \ref{delirio} and $f$ and $k$ are defined as above, we set
		\begin{align*}
			\Psi_x:=k_{\psi(x)}\circ(d\psi)_x\circ f_x^{-1}\colon\mathbb R^{2n+1}\to\mathbb R
			^{2n+1}\,,
		\end{align*}
		and we set
		\begin{align*}
			L_x:=k_{\psi(x)}^{-1}\circ \Psi_e\circ f_x\text{ and }R_x:=k_{\psi(x)}^{-1}\circ [\Psi_x-\Psi_e]\circ f_x\,.
		\end{align*}
		then 
		\begin{align*}
			(d\psi)_x=L_x+R_x\,
		\end{align*}
		where the map $L_x$ is an orthogonal transformation for any $x\in\mathbb H^{n}$, and the linear map
		$
			R_x
		$
		 is a smooth map 
		vanishing at $x=e$.
		Moreover
		\begin{equation}\label{chef}
			\ast\Lambda_h(d\psi)_x=\Lambda_{2n+1-h}(d\psi)_x\ast-\Lambda_{2n+1-h}R_x\ast+\ast\Lambda_hR_x\,.
		\end{equation}
		\end{lemma}
		\begin{proof}
		By Theorem \ref{delirio} $(d\psi)_e$ is an orthogonal linear map, thus
	\begin{align*}
		k_{\psi(e)}\circ(d\psi)_e\circ f_e^{-1}\colon \mathbb R
		^{2n+1}\to\mathbb R^{2n+1}\end{align*}
	is an orthonormal map, since by construction both $f_x$ and $k_m$ are orthonormal maps for any $x\in\mathbb H^n$ and for any $m\in M$, respectively, 
		that is $\Psi_e$ is an orthogonal linear transformation.
		
		Moreover, if we express
		\begin{align*}
			\Psi_x=\Psi_e+\Psi_x-\Psi_e
		\end{align*}
		we have
		\begin{align*}
			(d\psi)_x=k_{\psi(x)}^{-1}\circ \Psi_x\circ f_x=k_{\psi(x)}^{-1}\circ \Psi_e\circ f_x+k_{\psi(x)}^{-1}\circ [\Psi_x-\Psi_e]\circ f_x\,,
		\end{align*}
		so if we set
		\begin{align*}
			L_x:=k_{\psi(x)}^{-1}\circ \Psi_e\circ f_x\text{ and }R_x:=k_{\psi(x)}^{-1}\circ [\Psi_x-\Psi_e]\circ f_x\,.
		\end{align*}
		we obtain
		\begin{align*}
			(d\psi)_x=L_x+R_x\,\text{ and }\,
			\Lambda_h(d\psi)_x=\Lambda_hL_x+\Lambda_hR_x\,.
		\end{align*}

		In particular, we notice that 
		$
			R_x=(d\psi)_x-k_{\psi(x)}^{-1}\circ \Psi_e\circ f_x,
		$
		which can be seen as a matrix-valued smooth map 
		vanishing at $x=e$. 
		
		Moreover, the map $k_{\psi(x)}^{-1}\circ \Psi_e\circ f_x$ is an orthogonal transformation for any $x\in\mathbb H^{n}$, and by Proposition \ref{giaquinta}
		\begin{align*}
			\ast\Lambda_hL_x=\Lambda_{2n+1-h}L_x\ast\,.
		\end{align*}
		
		Therefore we can then write
		\begin{align*}\label{2pezzi}
			\ast\Lambda_h(d\psi)_x=&\ast\Lambda_hL_x+\ast\Lambda_hR_x=\Lambda_{2n+1-h}L_x\ast+\ast\Lambda_hR_x\\
			=&\Lambda_{2n+1-h}(d\psi)_x\ast-\Lambda_{2n+1-h}R_x\ast+\ast\Lambda_hR_x\,. 
		\end{align*}
	\end{proof}
	
	In the following lemma we discuss the interplay between the operations $\;\ast$ and pullback.
	\begin{lemma}\label{Lemma 1}
		Let $\alpha$ be a smooth form on $M$ of degree $h$, and let $\psi\colon B(e,1)\to M$ be as in Theorem \ref{delirio},
		then
		\begin{equation}\label{bellina}
			\psi^\sharp(\ast\alpha)-\ast\psi^\sharp\alpha=\sum_Ib_I\xi_I^\mathbb H\,,
		\end{equation}
		where $\{\xi_I^\mathbb H\}_I$ is a left-invariant  basis of the space of forms in $\mathbb H^n$ of degree $2n+1-h$, and the $b_I\in C^\infty(\mathbb H^n)$ are smooth coefficients that vanish when evaluated at the point $x=e$.
	\end{lemma}
	\begin{proof}
		Let $\alpha$ be a form on $M$ of degree $h$, and take $v_1\wedge\cdots\wedge v_{2n+1-h}\in\Lambda_{2n+1-h}\mathfrak{h}$ an arbitrary simple ${2n+1-h}$-vector of norm $\le 1$. If $x\in B(e,1)$, we can write, using \eqref{chef}  
		\begin{align*}
			&\langle(\psi^\sharp(\ast\alpha))_x\mid v_1\wedge\cdots\wedge v_{2n+1-h}\rangle=\langle \ast \alpha_{\psi(x)}\mid\Lambda_{2n+1-h}(d\psi)_x(v_1\wedge\cdots\wedge v_{2n+1-h})\rangle\\&
			=\langle\alpha_{\psi(x)}\mid\ast[\Lambda_{2n+1-h}(d\psi)_x(v_1\wedge\cdots\wedge v_{2n+1-h})]\rangle=\langle\alpha_{\psi(x)}\mid\Lambda_{h}(d\psi)_x\ast(v_1\wedge\cdots\wedge v_{2n+1-h})\rangle
			\\
			&-\langle\alpha_{\psi(x)}\mid\Lambda_{h}R_x\ast(v_1\wedge\cdots\wedge v_{2n+1-h})\rangle+\langle\alpha_{\psi(x)}\mid\ast\Lambda_{{2n+1-h}}R_x (v_1\wedge\cdots\wedge v_{2n+1-h})\rangle
			\\&=\langle (\ast\psi^\sharp\alpha)_x\mid (v_1\wedge\cdots\wedge v_{2n+1-h})\rangle-\langle\alpha_{\psi(x)}\mid\Lambda_{h}R_x\ast(v_1\wedge\cdots\wedge v_{2n+1-h})\rangle
			\\&+\langle\alpha_{\psi(x)}\mid\ast\Lambda_{{2n+1-h}}R_x (v_1\wedge\cdots\wedge v_{2n+1-h})\rangle\,.
		\end{align*}
		
		Hence
		\begin{equation*}
			\psi^\sharp (\ast \alpha) - \ast  \psi^\sharp  \alpha=: \sum_I b_I\xi_I^{\mathbb H},
		\end{equation*}
		where the $\xi_I^{\mathbb H}$  belong to a basis of the space of
		the forms of degree ${2n+1-h}$ in $\mathbb{H}^n$.
		
		Let $W_I^{\mathbb H}$ be the dual of the ${2n+1-h}\,$-covector $\xi_I^{\mathbb H}$. Then
		\begin{equation}\label{1}\begin{split}
				b_I (x) & =  \Scal{  (\psi^\sharp(\ast \alpha))_x - (\ast  \psi^\sharp  \alpha)_x}{W_I^{\mathbb{H}}} \\ &
				= \Scal{\alpha_{\psi(x)}}{ (\Lambda_{h}(R_x )\big(\ast W_I^{\mathbb{H}}\big)}
				-
				\Scal{ \alpha_{\psi(x)}}{\ast \Lambda_{{2n+1-h}}(R_x )W_I^{\mathbb{H}}}.
		\end{split}\end{equation}
	Since $ R$  vanishes at $x=e$ (see  Lemma \ref{RMK MAPPA}), then  $b_I(e)=0$.
		
		\end{proof}
		Using $\delta_c^M\alpha$ instead of $\alpha$ in \eqref{bellina} we get:
		{\begin{corollary}With the same hypotheses of Lemma \ref{Lemma 1} we have
		\begin{equation}\label{cavolata}
			\psi^\sharp(\ast\delta_c^M\alpha)-\ast\psi^\sharp(\delta_c^M\alpha)=:\sum_JB_J\xi_J^{\mathbb H}\,,
		\end{equation}
		where $\{\xi_J^\mathbb H\}_J$ is a left-invariant  basis of the space of forms in $\mathbb H^n$ of degree $2n+2-h$, and the $B_I\in C^\infty(\mathbb H^n)$ are smooth coefficients defined by
		\begin{equation}\label{capitalB}\begin{split}
				B_J (x) & =  \Scal{ (\psi^\sharp (\ast \delta_c^M \alpha))_x - (\ast  \psi^\sharp  (\delta_c^M\alpha))_x}{W_J^{\mathbb{H}}} \\ &
				= \Scal{(\delta_c^M\alpha)_{\psi(x)}}{ (\Lambda_{h-1}(R_x )\big(\ast W_J^{\mathbb{H}}\big)}
				-
				\Scal{ (\delta_c^M\alpha)_{\psi(x)}}{\ast \Lambda_{2n+2-h}(R_x )W_J^{\mathbb{H}}}.
		\end{split}\end{equation}
		and 		
		that vanish when evaluated at the point $x=e$.
		\end{corollary}}
		
		\begin{remark}\label{12lugliobis}
 Let us denote by $r_{j,x}$ the coefficients of $\Lambda_{h}(R_x )W_I^{\mathbb{H}}$. The  $r_{j,x}$ are smooth maps that vanish at $x=e$.  
		Let us notice that, even though the coefficients $b_I$ are functions on $\mathbb H^n$, by the second equality of \eqref{1} in the proof above, they can be expressed as a linear combination of terms of the form
		\begin{align}\label{espressione di b}
			\alpha_{i,\psi(x)}\,r_{j,x}=(\alpha_{i}\circ\psi)(x)\,r_{j}(x)
		\end{align}
			where we used the subscript  $\psi$ to highlight the dependence of $b_I$ on the map $\psi$.
			\end{remark}

		\begin{remark}\label{13luglio}
Let us assume $\alpha$ to be a form on $M$ of degree $h$, then by Lemma \ref{Lemma 1} we know
		\begin{align}\label{bellina3}
			d_c^\mathbb H\psi^\sharp(\ast\alpha)-d_c^\mathbb H\ast\psi^\sharp(\alpha)=d_c^\mathbb H\sum_Ib_I\xi_I^\mathbb H\,.
		\end{align}
		Applying the Hodge operator then we get
		$$
			\ast d_c^\mathbb H\psi^\sharp(\ast\alpha)-\ast d_c^\mathbb H\ast\psi^\sharp(\alpha)=\ast d_c^\mathbb H\sum_Ib_I\xi_I^\mathbb H\,.
		$$
		Keeping in mind that  $d_c^\mathbb H\psi^\sharp=\psi^\sharp d_c^M$ the expression above becomes
		$
			\ast\psi^\sharp(d_c^M\ast\alpha)-\delta_c^\mathbb H\psi^\sharp(\alpha)=\ast d_c^\mathbb H\sum_Ib_I\xi_I^\mathbb H\,.
		$
		Therefore, writing, up to a sign, $1=\ast\ast$, from the last equality we have $
			\ast\psi^\sharp(\ast\ast d_c^M\ast\alpha)-\delta_c^\mathbb H\psi^\sharp(\alpha)=\ast d_c^\mathbb H\sum_Ib_I\xi_I^\mathbb H\,,
		$
		and, again up to a sign,  we eventually get
		\begin{align}\label{bellina2}
			\delta_c^\mathbb H\psi^\sharp(\alpha)=\ast\psi^\sharp(\ast\delta_c^M\alpha)+\ast d_c^\mathbb H\sum_Ib_I\xi_I^\mathbb H\,.
		\end{align}
		
	If $2n+1-h\neq n$, i.e $h\neq n+1$ , the differential $d_c^\mathbb H$ has order $1$, then by left-invariance, we have 
		$$
		\ast d_c^{\mathbb H} \sum_Ib_I\xi_I^\mathbb H=\sum_{\ell, I} (W_\ell^{\mathbb H} b_I) \ast(\omega^\mathbb H_\ell\wedge \xi_I)\,.
		$$
		Then $\ast d_c^\mathbb H\sum_Ib_I\xi_I^\mathbb H$ is a form of degree $h-1$ whose coefficients are of the type
		\begin{align*}
		 	(W_\ell^{\mathbb H} \alpha_{i,\psi(x)})\cdot r_{j,x} \;\text{ or }\;	\alpha_{i,\psi(x)}\cdot (W_\ell ^{\mathbb H}r_{j,x})\;\;\,.
		\end{align*}
		When {$h= n+1$} the differential $d_c^{\mathbb H}$ has order two. Then
		 $\ast d_c^\mathbb H\sum_Ib_I\xi_I^\mathbb H$ is a form of degree $n$ whose coefficients are of type 
		\begin{align*}
			(W_\ell^{\mathbb H}W_\lambda^{\mathbb H} \alpha_{i,\psi(x)})\cdot r_{j,x} \;,\;(W_\ell^{\mathbb H}\alpha_{i,\psi(x)})\cdot W_\lambda^{\mathbb H} r_{j,x}\;\text{ or }\;\alpha_{i,\psi(x)}\cdot (W_\ell^{\mathbb H} W_\lambda^{\mathbb H} r_{j,x})\;\,.
		\end{align*}
		
		{Moreover, by \eqref{capitalB}, the coefficients $B_J$  are of the type
		\begin{align*}
		 	(W_\ell^{\mathbb H} \alpha_{i,\psi(x)})\cdot r_{j,x} \;\text{ if }\; h\neq n+1,	\text{ or }\ (W_\ell^{\mathbb H}W_\lambda^{\mathbb H} \alpha_{i,\psi(x)})\cdot r_{j,x}\;\text{ if  }\;h=n+1\, .\end{align*}}
\end{remark}


\bigskip

We are now in position to examine the  interplay between $\psi^\sharp\delta_c^M $ and $\delta_c^{\mathbb H}\psi^\sharp $.

{
\begin{proposition}\label{interplay} With the same hypotheses of Lemma \ref{Lemma 1} we have
\begin{equation}\label{quella finale}
	\delta_c^\mathbb H\psi^\sharp(\alpha)=\psi^\sharp(\delta_c^M\alpha)+\ast d_c^\mathbb H\sum_Ib_I\xi_I^\mathbb H\,+\ast\sum_JB_J\xi_J^\mathbb H.
\end{equation}
where the coefficient $b_I$ and $B_J$ are defined in \eqref{1} and \eqref{capitalB} respectively.
\end{proposition}
\begin{proof} We start from \eqref{bellina2} and we combine with \eqref{cavolata} to get, , up to a sign,
$$\delta_c^\mathbb H\psi^\sharp(\alpha)=\ast\ast\psi^\sharp(\delta_c^M\alpha)+\ast d_c^\mathbb H\sum_Ib_I\xi_I^\mathbb H+\ast\sum_JB_J\xi_J^\mathbb H\,.$$
\end{proof}}

Given a smooth $h-$form $\alpha$ on $M$, we  show below that the $L^p$  norm of $$\delta_c^{\mathbb H}\psi^\sharp \alpha-\psi^\sharp\delta_c^M\alpha$$  is small if  $\alpha$ is  supported in a  suitably chosen small ball.
 In order to show this, we have to handle terms of the form $(W_\ell^{\mathbb H} (\alpha_{i,\psi(x)})\cdot r_{j,x}$ and  terms of the form $\alpha_{i,\psi(x)}\cdot (W_\ell ^{\mathbb H}r_{j,x})$ (and with a little difference in the case $h=n+1$). The proof relies on two different approaches: for terms like $\alpha_{i,\psi(x)}\cdot (W_\ell ^{\mathbb H}r_{j,x})$ we  a Sobolev inequality  and 
the fact that the support of $ \alpha $ is small; on the contrary when we want to handle terms of the form $(W_\ell^{\mathbb H} (\alpha_{i,\psi(x)})\cdot r_{j,x}$, we  use the fact that the  $ r_{j,x} $ tend to zero as $x\to e$. Remember that the $r_{j,x}$ are the coefficients of  $\Lambda_{h}(R_x )W_I^{\mathbb{H}}$, which depend only on the map $\psi$ (see Lemma \ref{RMK MAPPA}), and the  map $\psi$ can be controlled with constants depending only on the geometry of $M$ i.e. on $r$ and $C_M$ (by Theorem \ref{delirio}). In conclusion we find that the radius of the support can be chosen independently of $\alpha$ and depends only on  the geometry of $M$.

\bigskip

\begin{proposition}\label{Proposition P}  For any $a\in M$ we consider the map $\psi=\psi_a$ as in Theorem  \ref{delirio}.
With the notation of Proposition \ref{interplay} we define the operator
		\begin{align}\label{superP}
			\mathcal{P}(x,W^{\mathbb H})\psi^\sharp\alpha:=\delta_c^\mathbb H(\psi^\sharp\alpha)-\psi^\sharp(\delta_c^M\alpha)=\ast d_c^\mathbb H\sum_Ib_I\xi_I^\mathbb H+\ast\sum_JB_J\xi_J^\mathbb H\ \,.
		\end{align}
		The operator $\mathcal{P}$
			 is  a linear differential operator on $\he n$ which
			is of second order if $h= n+1$, and of first order otherwise.

Let $1<p<\infty$ and let  $\varepsilon >0$. Then there exists $\overline\eta=\overline\eta(C_M, r, \varepsilon)>0$ such that if $\eta<\overline\eta(C_M, \varepsilon)$ and  $\alpha$ is a smooth  $h$-form on $M$ supported in $B(a,\eta)$ 
we have
				\begin{align}\label{vigilia}
			\|\mathcal{P}(\psi^\sharp\alpha)\|_{L^p( \psi^{-1}(B(a,\eta), E_0^{h-1})}\le \varepsilon  \|\psi^\sharp\alpha\|_{W^{1,p}(\psi^{-1}(B(a,\eta)), E_0^h)}\,,
		\end{align}
		if $h\neq n+1$, or
		\begin{align}\label{vigilia2}
			\|\mathcal{P}(\psi^\sharp\alpha)\|_{L^p( \psi^{-1}(B(a,\eta)), E_0^{n})}\le \varepsilon  \|\psi^\sharp\alpha\|_{W^{2,p}(\psi^{-1}(B(a,\eta)), E_0^{n+1})}\,.
		\end{align}

	\end{proposition}

\begin{proof} {By Proposition \ref{interplay},}  the operator $\mathcal P$ is an operator acting on $\psi^\sharp\alpha$, since  we can write
{\begin{align*}
			\mathcal{P}(x,W^{\mathbb H})\psi^\sharp\alpha=\delta_c^\mathbb H(\psi^\sharp\alpha)-\psi^\sharp\delta_c^M(\psi^\sharp)^{-1}(\psi^\sharp\alpha)\,.
		\end{align*}}

\medskip

		 As pointed out in Lemma \ref{Lemma 1}  $r_{j,x}$ are smooth maps that vanish at $x=e$ and  since they are the coefficients of $\Lambda_{2n+1-h}(R_x )W_I^{\mathbb{H}}$, they do not depend on $\alpha$ but only on $\psi$. 
		
 Keeping  into account \eqref{espressione di b} that expresses $b_I$, we get estimates of the type
	$$
|b_I(x)|\le |\alpha_{i,\psi(x)}| \Omega(x)
$$
where $\Omega(x)=O(\rho(x))$ for $x\to e$ (with a slight abuse of notation here and in the sequel we avoid to take the sum over the index $i$). 

Moreover if $h\neq n+1$,  from \eqref{bellina2}, by the triangular inequality, we also  get estimates of the type


{\begin{equation}\label{unouno}
|\delta_c^\mathbb H\psi^\sharp(\alpha)-\psi^\sharp(\delta_c^M\alpha)|\le \sum_\ell|W_\ell^{\mathbb{H}}\alpha_{i,\psi(x)}|\Omega(x)+|\alpha_{i,\psi(x)}\cdot \sum_\ell W_\ell ^{\mathbb H}r_{j,x}|\,,
 \end{equation}
where $\Omega(x)=O(\rho(x))$ for $x\to e$.}

Similarly if  $h= n+1$, remembering also that $d_c^M$ is a second order differential operator, we have estimates of the type

\begin{equation}\label{duedue}\begin{split}
|\delta_c^\mathbb H&\psi^\sharp(\alpha)-\psi^\sharp(\delta_c^M\alpha)|\le \sum_{\ell,\lambda}|W_\ell^{\mathbb{H}}W_\lambda^{\mathbb{H}}\alpha_{i,\psi(x)}|\Omega(x)
\\&
+ \sum_{\ell}|W_\ell^{\mathbb H}\alpha_{i,\psi(x)}\cdot\sum_\lambda W_\lambda^{\mathbb H} r_{j,x} |+|\alpha_{i,\psi(x)}\cdot \sum_{\ell,\lambda}(W_\ell^{\mathbb H} W_\lambda^{\mathbb H}) r_{j,x}| \,,
\end{split}
\end{equation}
where $\Omega(x)=O(\rho(x))$ for $x\to e$.

\begin{itemize}
	\item Case $h\neq n+1$.
	
Let us estimate the $L^p$ norm of the first term in the right hand side of \eqref{unouno}.
If $x\to e$ then $\Omega(x)\to 0$,   then the $L^p$-norm of term $\sum_\ell|W_\ell^{\mathbb H}\alpha_{i,\psi(x)}|\Omega(x)$ is controlled by $\varepsilon \sum_\ell\|W_\ell^{\mathbb H}\psi^\sharp\alpha\|_{L^p}$ provided $x$ is sufficiently close to $e$. Thus, now we need to  estimate  the second term of the right hand side of \eqref{unouno}.

		Hence, we need only to handle carefully the terms that can be expressed as a linear combination of terms the form
		\begin{align*}
				\alpha_{i,\psi(x)}\cdot (W_\ell ^{\mathbb H}r_{j,x})=\alpha_{i}\circ\psi(x)\cdot (W_\ell ^{\mathbb H}r_{j,x})
		\end{align*}
		
		The functions 
		${W_\ell ^{\mathbb H}r_{j,x}}
		$
		are bounded in $\psi^{-1}(B(a,\eta)) $. But $B(e,\frac{\eta}{C_M})\subseteq\psi^{-1}(B(a,\eta))\subseteq B(e,C_M{\eta})$
		and by the Sobolev inequality 
		$$
		\left(\int_{B(e,C_M\eta)}(\alpha_{i}\circ \psi(x))^p dx\right)^{1/p}\le c_pC_M\eta \left(\int_{B(e,C_M\eta)}\sum_\ell|W^{\mathbb H}_\ell(\alpha_{i}\circ\psi)(x)|^p dx\right)^{1/p}
		$$
		where $c_p$ denotes the Sobolev constant (depending only on $p$ and $n$).
		
		If we chose $\eta$ so that $c_pC_M\eta<\varepsilon$ i.e. $\eta<\overline\eta$
		\begin{equation}\label{eta}
		\overline\eta=\frac{\varepsilon}{c_pC_M}
		\end{equation}
		finally we get,  
		$$
		\|\alpha_{i,\psi(x)}\cdot (W_\ell ^{\mathbb H}r_{j,x})\|_{L^p(\psi^{-1}(B(a,\eta)))}
		\le \varepsilon \sum_\ell\|W^{\mathbb H}_\ell\alpha_{i,\psi}\|_{L^p(\psi^{-1}(B(a,\eta)))}\,.
		$$
		Therefore reasoning on the differential form $\alpha$, and possibly relabeling $\varepsilon$,  we get \eqref{vigilia}.

		\item Case $h=n+1$. Arguing as above, 
						again,  we notice that the first term on the right hand side of \eqref{duedue} is
		 of the form $(W_\ell^{\mathbb H}W_\lambda^{\mathbb H} (\alpha_{i,\psi(x)})\cdot r_{j,x}$ and can be estimated 
		by $\varepsilon \sum_{\ell,\lambda}\|W_\ell^{\mathbb H}W_\lambda^{\mathbb H}(\psi^\sharp\alpha)\|_{L^p}$ since $r_{j,x}\to 0$ if $x\to e$.
		
		Hence, we need only to handle carefully the terms that can be expressed as a linear combination of terms the form
		\begin{align*}
			 \;\;(W_\ell^{\mathbb H}\alpha_{i,\psi(x)})\cdot W_\lambda^{\mathbb H}) r_{j,x}\;\text{ and }\;\alpha_{i,\psi(x)}\cdot (W_\ell^{\mathbb H} W_\lambda^{\mathbb H}) r_{j,x}\;.
		\end{align*}
	They	can be handled again by using the Sobolev inequality, since both first and second derivatives of $r_{j,x}$ are bounded in $B(e,C_M\eta)$. Indeed, when we apply Sobolev inequality  to terms of the forms $W_\ell^{\mathbb H}\alpha_{i,\psi(x)}$  we get an estimate of the form  $c_pC_M\eta \sum_{\ell,\lambda}\|W_\ell^{\mathbb H}W_\lambda^{\mathbb H}(\psi^\sharp\alpha)\|_{L^p( \psi^{-1}(B(a,\eta)))}$. Similarly, the terms of the type $\alpha_{i,\psi(x)}\cdot (W_\ell^{\mathbb H} W_\lambda^{\mathbb H}) r_{j,x}$ can be estimated by $c_pC_M\eta\sum_\ell\|W^{\mathbb H}_\ell(\psi^\sharp\alpha)\|_{L^p( \psi^{-1}(B(a,\eta)))}$.
			
			Therefore, again choosing $c_pC_M\eta\le \varepsilon$, eventually we get \eqref{vigilia2}
			
		\end{itemize}
	\end{proof}
	Notice that  $\overline \eta$ depends also on $p, n$, but it is not explicit  in the statement above since it is well known and what is relevant to us is to show  the dependence on the geometry of $M$.
	\begin{remark}\label{il caso h=n}
	If the differential form $\alpha$ is of degree $n$, in order to prove Theorem \ref{gaffney-M},  we shall also need to know the interplay between $\psi^\sharp(d_c^M\delta_c^M) $ 
	and $d_c^{\mathbb H}\delta_c^{\mathbb H}\psi^\sharp $. 
	
	We set
	\begin{align}\label{superQ}
	\mathcal{Q}(x,W^{\mathbb H})(\psi^\sharp\alpha):=d_c^{\mathbb H}\delta_c^{\mathbb H}(\psi^\sharp\alpha)-\psi^\sharp(d_c^M\delta_c^M\alpha)\,.
	\end{align}
	Now, $\mathcal{Q}(\psi^\sharp\alpha):=d_c^{\mathbb H}\delta_c^{\mathbb H}(\psi^\sharp\alpha)-\psi^\sharp(d_c^M\delta_c^M\alpha)=d_c^{\mathbb H}\delta_c^{\mathbb H}(\psi^\sharp\alpha)-d_c^{\mathbb H}\psi^\sharp(\delta_c^M\alpha)$. Hence,   keeping in mind  \eqref{superP}, we easily see that $$\mathcal{Q}(x,W^{\mathbb H})=d_c^\mathbb H\mc P(x,W^{\mathbb H})$$  is a second order differential operator (since $\mc P\psi^\sharp \alpha $ is a form of degree $n-1$ and hence $d^\mathbb H$ is a differential operator of degree $1$). To estimate the $L^p $ norm of $\mathcal{Q}(x,W^{\mathbb H})(\psi^\sharp\alpha)$ we have to estimates terms of  type 
		\begin{align*}
			(W_\ell^{\mathbb H}W_\lambda^{\mathbb H} (\alpha_{i,\psi(x)}))\cdot r_{j,x} \;,\;(W_\ell^{\mathbb H}\alpha_{i,\psi(x)})\cdot W_\lambda^{\mathbb H} r_{j,x}\;\text{ or }\;\alpha_{i,\psi(x)}\cdot (W_\ell^{\mathbb H} W_\lambda^{\mathbb H}) r_{j,x}\;\,.
		\end{align*}
	Hence, with the notation of Proposition \ref{Proposition P}, given $\varepsilon >0$ if $\eta<\overline\eta(C_M, \varepsilon)$ and
	  $\alpha$ is supported in $B(a,\eta)$  then
		
		\begin{align}\label{Q}
			\|\mathcal{Q}(\psi^\sharp\alpha)\|_{L^p}\le \varepsilon  \|\psi^\sharp\alpha\|_{W^{2,p}}\,.
		\end{align}
		
		Likewise, if the differential form $\alpha$ is of degree $n+1$, in proving  Theorem \ref{gaffney-M},  we shall need also to evaluate the difference $\delta_c^{\mathbb H}d_c^{\mathbb H}(\psi^\sharp\alpha)-\psi^\sharp(\delta_c^Md_c^M\alpha)\,.$
		Let us set 
		$$
	\mathcal{T}(x,W^{\mathbb H})(\psi^\sharp\alpha):=\delta_c^{\mathbb H}d_c^{\mathbb H}(\psi^\sharp\alpha)-\psi^\sharp(\delta_c^Md_c^M\alpha)\,,
	$$
	hence $$\mathcal{T}(x,W^{\mathbb H})=\mc P(x,W^{\mathbb H}) d_c^\mathbb H$$
	 is again a second order differential operator (since $d_c^\mathbb H\psi^\sharp \alpha $ is a form of degree $n+2$ and hence $\mc P$ is a differential operator of degree $1$). Analogously as the operator\eqref{Q}, we have
	\begin{align}\label{T}
			\|\mathcal{T}(\psi^\sharp\alpha)\|_{L^p}\le \varepsilon  \|\psi^\sharp\alpha\|_{W^{2,p}}\,.
		\end{align}
		
		The estimates obtained in Proposition \ref{Proposition P}, the estimates \eqref{Q} and \eqref{T} will be used in the proof of Theorem \ref{lemma ivaniec}.
	\end{remark}

	In the sequel, in order to prove a Gaffney estimate on $M$, once we have chosen an atlas on $M$ and a partition of the unity subordinated to the chosen covering of $M$,  we will need also to obtain $L^p$ estimates of the commutator between the operator $d_c^M$ (or $\delta_c^M$) with a smooth function. To this aim it will be useful the following lemma (since we are interested in obtaining  $L^p$ estimates,  the following equalities  are true up to signs).
	\begin{lemma}\label{commutator bis}
	Let $\psi $ be a contactomorphism  from an open set $\mc U\subset \he{n}$ to $M$, and we denote by $\mc V$
the open set  $\mc V= \phi(\mc U)$. 	
	If $\chi$ is a smooth function in $M$ and $\alpha\in E_0^h(\mc V)$, we have 
	\begin{align}
		\begin{split}\label{pullback commutator i.}
			\psi^\sharp\big([d_c^M,\chi]\alpha\big)=&[d_c^\mathbb H,\chi\circ\psi]\psi^\sharp\alpha\,, \quad\quad \mathrm{for \ any }\ h
		\end{split}\\
		\begin{split}\label{pullback commutator ii.}
			{\psi^\sharp\big([\delta_c^M,\chi]\alpha\big)}=&{[\delta_c^\mathbb H,  \chi\circ \psi]\psi^\sharp\alpha+[\mc P,\chi\circ\psi]\psi^\sharp(\ast\alpha)}\,, \quad\quad \mathrm{if}\ h\neq n+1
		\end{split}\\
		\begin{split}\label{pullback commutator iii.}
			\psi^\sharp\big([d_c^M\delta_c^M,\chi]\alpha\big)=&[d_c^\mathbb H\delta_c^\mathbb H \chi\circ\psi]\psi^\sharp\alpha+[d_c^\mathbb H\mc P,\chi\circ\psi]\psi^\sharp\alpha\,, \quad\quad \mathrm{if}\ h= n
		\end{split}\\
		\begin{split}\label{pullback commutator iv.}
			\psi^\sharp\big([\delta_c^M d_c^M,\chi]\alpha\big)=&[\delta_c^\mathbb H d_c^\mathbb H,\chi\circ\psi]\psi^\sharp\alpha+[\mc Pd_c^\mathbb H,\chi\circ\psi^\sharp]\psi^\sharp\alpha\,\quad\quad \mathrm{if}\ h= n+1.
		\end{split}
	\end{align}
%
\end{lemma}
\begin{proof}
	The first equality follows directly from the fact that the Rumin differential and the pullback of a contact map commute (see also Proposition \ref{memory}):
	\begin{align*}
		\psi^\sharp\big([d_c^M,\chi]\alpha\big)=&\psi^\sharp d_c^M(\chi\alpha)-\psi^\sharp\big(\chi\ d_c^M\alpha\big)=d_c^\mathbb H\psi^\sharp(\chi\alpha)-\chi\circ\psi\cdot\psi^\sharp (d_c^M \alpha)=\\=&d_c^\mathbb H(\chi\circ\psi\cdot\psi^\sharp\alpha)-\chi\circ\psi\cdot d_c^\mathbb H\psi^\sharp\alpha=[d_c^\mathbb H,\chi\circ \psi]\psi^\sharp\alpha\,.
	\end{align*}
	For the second formula, the codifferential and the pullback map do not commute, however we can use   \eqref{superP}:
	
	{\begin{align*}
		\psi^\sharp\big(&[\delta_c^M,\chi]\alpha\big)=\psi^\sharp\delta_c^M(\chi\alpha)-\psi^\sharp(\chi\delta_c^M\alpha)
		\\
		=&\delta_c^\mathbb H\psi^\sharp(\chi\alpha)+\mc P\psi^\sharp(\chi\alpha)-\chi\circ\psi\big(\psi^\sharp\delta_c^M\alpha\big)
		\\
		=&\delta_c^\mathbb H\big(\chi\circ\psi\cdot \psi^\sharp\alpha\big)+\mc P(\chi\circ\psi\cdot\psi^\sharp\alpha)-\chi\circ\psi\big(\delta_c^\mathbb H\psi^\sharp\alpha+\mc P\psi^\sharp\alpha\big)\\=&[\delta_c^\mathbb H,\chi_j\circ\psi_j]\psi_j^\sharp\alpha+[\mc P,\chi\circ\psi]\psi^\sharp\alpha\,.
	\end{align*}
	
	}

	The third and fourth formulae will follow by using a similar reasoning as before:
	\begin{align*}
		\psi^\sharp&\big([d_c^M\delta_c^M,\chi]\alpha\big)=\psi^\sharp\big(d_c^M\delta_c^M(\chi\alpha)\big)-\psi^\sharp\big(\chi d_c^M\delta_c^M\alpha\big)\\=&d_c^\mathbb H\psi^\sharp(\delta_c^M(\chi\alpha))-\chi\circ\psi\  d_c^\mathbb H\psi^\sharp(\delta_c^M\alpha)
		\\
		=&d_c^\mathbb H\big(\delta_c^\mathbb H\psi^\sharp(\chi\alpha)+\mc P\psi^\sharp(\chi\alpha)\big)-\chi\circ\psi\ d_c^\mathbb H\big(\delta_c^\mathbb H\psi^\sharp\alpha+\mc P\psi^\sharp\alpha\big)
		\\
		=&d_c^\mathbb H\big(\delta_c^\mathbb H\big(\chi\circ\psi\cdot \psi^\sharp\alpha\big)+\mc P(\chi\circ\psi\cdot\psi^\sharp\alpha)\big)-\chi\circ\psi\  d_c^\mathbb H\big(\delta_c^\mathbb H\psi^\sharp\alpha+\mc P\psi^\sharp\alpha\big)
		\\=&[d_c^\mathbb H\delta_c^\mathbb H,\chi\circ\psi]\psi^\sharp\alpha+[d_c^\mathbb H\mc P,\chi\circ\psi]\psi^\sharp\alpha
	\end{align*}
	
	
	and as for \eqref{pullback commutator iv.}, using again \eqref{superP}, 
	\begin{align*}
		\psi^\sharp&\big([\delta_c^M d_c^M,\chi]\alpha\big)=\psi^\sharp\big(\delta_c^Md_c^M(\chi\alpha)\big)-\psi^\sharp\big(\chi\delta_c^M d_c^M\alpha\big)\\=&\delta_c^\mathbb H\psi^\sharp d_c^M(\chi\alpha)+ \mc P\psi^\sharp d_c^M(\chi\alpha)-\chi\circ\psi\big(\delta_c^\mathbb H\psi^\sharp d_c^M\alpha+\mc P\psi^\sharp d_c^M\alpha\big)
		\\=&
		\delta_c^\mathbb H d_c^\mathbb H\psi^\sharp(\chi\alpha)+ \mc Pd_c^\mathbb H\psi^\sharp (\chi\alpha)-\chi\circ\psi\big(\delta_c^\mathbb Hd_c^\mathbb H\psi^\sharp \alpha+\mc P d_c^\mathbb H\psi^\sharp\alpha\big)
		\\=&[\delta_c^\mathbb H d_c^\mathbb H,\chi\circ\psi]\psi^\sharp\alpha+ [\mc Pd_c^\mathbb H,\chi\circ \psi]\psi^\sharp\alpha
	\end{align*}

\end{proof}

\section{Sobolev-Gaffney type inequalities on contact manifolds}

{
We recall the following  Sobolev-Gaffney type inequalities proved in the setting of Heisenberg groups for differential forms in $\mc D(\he{n}, E_0^h)$ (see \cite{BF7}, Remark 5.3, i), iii),vi) therein) which by density of $\mc D(\he{n}, E_0^h)$ in $W^{k,p}(\he{n}, E_0^h)$ (see Theorem \ref {denso in h}) can be rephrased as:
\begin{lemma}\label{gaffney-Hn}
Let $1\le h\le 2n$,   and $1<p<\infty$, then there
exists a constant $C_G=C_G(p,n,h)>0$ such that for all $u \in {W^{1,p}(\he{n}, E_0^h)}$ we have
\begin{itemize}

\item[i)]  
\begin{equation}\label{p>1:eq1}\begin{split}
\| u \|_{W^{1,p}(\he{n}, E_0^h)} &
\le C_G\big( \| d_c^{\mathbb{H}} u \|_{L^{p}(\he{n}, E_0^{h+1})} +
\| \delta_c^{\mathbb{H}} u \|_{L^{p}(\he{n}, E_0^{h-1})}
\\& \hphantom{xxxxx}+ \|u\|_{L^p(\he{n}, E_0^h)}
\big),
\end{split}\end{equation}
if $h\neq n, n+1$;

\medskip

\item[ii)]  if $h= n$,
\begin{equation}\label{treno}\begin{split}
\| u \|_{W^{2,p}(\he{n}, E_0^n)} &
\le C_G\big( \| d_c^{\mathbb{H}}u \|_{L^{p}(\he{n}, E_0^{n+1})} +
\| d_c^{\mathbb{H}}\delta_c^{\mathbb{H}} u \|_{L^{p}(\he{n}, E_0^{n})}
\\&
 \hphantom{xxxxx} +  \|u\|_{L^p(\he{n}, E_0^n)}\big),
\end{split}\end{equation}

 %
%
%
%
%
and,
\item[iii)] if $h= n+1$, 
\begin{equation}\label{p>1:eq2bis}\begin{split}
\| u \|_{W^{2,p}(\he{n}, E_0^{n+1})} &
\le C_G\big( \| \delta_c^{\mathbb{H}}d_c^{\mathbb{H}}u \|_{L^{p}(\he{n}, E_0^{n+1})} +
\| \delta_c^{\mathbb{H}} u \|_{L^{p}(\he{n}, E_0^{n})}
\\&
 \hphantom{xxxxx} +  \|u\|_{L^p(\he{n}, E_0^{n+1})}\big),
\end{split}\end{equation}
\end{itemize}

 %
%

\end{lemma}
}

Before stating the global result, we prove the following local one where we can use
the groundwork  
 just developed, together with the Gaffney-Sobolev inequality stated in the previous result.

\begin{theorem}\label{lemma ivaniec}
 Let $1<p<\infty$, there exists a positive constant $\tilde\eta=\tilde\eta(C_M,C_G)$ so that, if  $\eta<\tilde\eta$  and $(B(a_\ell,\eta), \psi_{a_\ell})$ is a chart of the atlas given in Remark \ref{covering ad hoc}, and $\alpha$ is a smooth form in $M$ with  support contained in $B(a_\ell,\eta)$, then there exists a constant $C=C(C_M, C_G)$  so that, if $h\neq n, n+1$
\begin{align}\label{local ineq to prove}
\begin{split}
			\Vert \psi_{a_\ell}^\sharp\alpha &\Vert_{W^{1,p}(\psi_{a_\ell}^{-1}(B(a_\ell,\eta)),E_0^h)}\le C\big(\Vert \psi_{a_\ell}^\sharp\alpha\Vert_{L^p(\psi_{a_\ell}^{-1}(B(a_\ell,\eta)),E_0^h)}
			\\
			&
+\Vert \psi_{a_\ell}^\sharp(d_c^M\alpha)\Vert_{L^p(\psi_{a_\ell}^{-1}(B(a_\ell,\eta)),E_0^{h+1})}+\Vert\psi_{a_\ell}^\sharp(\delta_c^M\alpha)\Vert_{L^p(\psi_{a_\ell}^{-1}(B(a_\ell,\eta)),E_0^{h-1})}\big)\,.
\end{split}
		\end{align}
		Whereas, if $h= n$ we get
		\begin{align}\label{local ineq to prove due}
		\begin{split}
			\Vert \psi_{a_\ell}^\sharp\alpha &\Vert_{W^{2,p}(\psi_{a_\ell}^{-1}(B(a_\ell,\eta)),E_0^n)}\le C\big(\Vert \psi_{a_\ell}^\sharp\alpha\Vert_{L^p(\psi_{a_\ell}^{-1}(B(a_\ell,\eta)),E_0^n)}
			\\
			&
+\Vert \psi_{a_\ell}^\sharp(d_c^M\alpha)\Vert_{L^p(\psi_{a_\ell}^{-1}(B(a_\ell,\eta)),E_0^{n+1})}+\Vert\psi_{a_\ell}^\sharp(d_c^M\delta_c^M\alpha)\Vert_{L^p(\psi_{a_\ell}^{-1}(B(a_\ell,\eta)),E_0^{n})}\big)\,.
\end{split}
		\end{align}
		and if  $h= n+1$ we get
		\begin{align}\label{local ineq to prove tre}
		\begin{split}
			\Vert \psi_{a_\ell}^\sharp\alpha &\Vert_{W^{2,p}(\psi_{a_\ell}^{-1}(B(a_\ell,\eta)),E_0^{n+1})}\le C\big(\Vert \psi_{a_\ell}^\sharp\alpha\Vert_{L^p(\psi_{a_\ell}^{-1}(B(a_\ell,\eta)),E_0^{n+1})}
			\\
			&
+\Vert \psi_{a_\ell}^\sharp(\delta_c^Md_c^M\alpha)\Vert_{L^p(\psi_{a_\ell}^{-1}(B(a_\ell,\eta)),E_0^{n+1})}+\Vert\psi_{a_\ell}^\sharp(\delta_c^M\alpha)\Vert_{L^p(\psi_{a_\ell}^{-1}(B(a_\ell,\eta)),E_0^{n})}\big)\,.
\end{split}
		\end{align}
\end{theorem}

\begin{proof}
We consider a chart $(B(a_\ell,\eta), \psi_{a_\ell})$ and we can assume that $\psi_{a_\ell}$ is taken as  in Theorem \ref{delirio}.  

To avoid cumbersome notation, in the sequel we  omit the subscripts and we  write $B(a,\eta)$ and $\psi$. Moreover, we set
$$\tilde B_\eta:=\psi^{-1}(B(a,\eta))$$
 and  write $L^p(\tilde B_\eta)$ instead of $L^p(\psi^{-1}(B(a,\eta)),E_0^{h})$ and similarly for the notation on Sobolev spaces.
If $\alpha$ is supported in $B(a,\eta)$ then,  
without loss of generality, we may assume that
$\psi^\sharp \alpha$ 
is compactly supported in $B(e,\eta/C_M)$ since $B(e,\frac{\eta}{C_M})\subseteq\psi^{-1}(B(a,\eta))\subseteq B(e,C_M{\eta})$.

		
		\medskip
		
		In order to prove \eqref{local ineq to prove}, we use first \eqref{p>1:eq1} with \eqref{superP} and then \eqref{vigilia}: given $\varepsilon>0$, by Proposition \ref{Proposition P} there exists $\overline\eta$ (see \eqref{eta}) so that, if $\eta<\overline\eta$ we obtain
		\begin{align*}
			\Vert&\psi^\sharp\alpha\Vert_{W^{1,p}(\tilde B_\eta)}\le C_G\big(\Vert\psi^\sharp\alpha\Vert_{L^p(\mathbb H^n)}+\Vert d_c^\mathbb H(\psi^\sharp\alpha)\Vert_{L^p(\mathbb H^n)}+\Vert\delta_c^\mathbb H(\psi^\sharp\alpha)\Vert_{L^p(\mathbb H^n)}\big)
			\\\le &C_G\big(\Vert\psi^\sharp\alpha\Vert_{L^p(\tilde B_\eta)}+\Vert \psi^\sharp(d_c^M\alpha)\Vert_{L^p(\tilde B_\eta)}+\Vert\psi^\sharp(\delta_c^M\alpha)\Vert_{L^p(\tilde B_\eta)}+\Vert\mathcal{P}(\psi^\sharp\alpha)\Vert_{L^p(\tilde B_\eta)}\big)
			\\
			\le &C_G\big(\Vert\psi^\sharp\alpha\Vert_{L^p(\tilde B_\eta)}+\Vert \psi^\sharp(d_c^M\alpha)\Vert_{L^p(\tilde B_\eta)}+\Vert\psi^\sharp(\delta_c^M\alpha)\Vert_{L^p(\tilde B_\eta)}+\varepsilon \Vert\psi^\sharp\alpha\Vert_{W^{1,p}(\tilde B_\eta)}\big)\,.
		\end{align*}
Choosing $\varepsilon\le 1/(2C_G)$ we have proved \eqref{local ineq to prove} for $\eta<\tilde\eta:=\frac{1}{2c_pC_MC_G}$ (notice that in the statement, again,  the dependence on $c_p$ was omitted). 

\bigskip

Let now $h=n$. The argument above needs to be only slightly modified. Indeed, we will apply the Gaffney inequality \eqref{treno},  where both $d_c^{\mathbb H}$ and $d_c^{\mathbb H}\delta_c^{\mathbb H} $ appearing on the right hand side are differential operators of order $2$. 
Therefore,
keeping also in mind \eqref{Q}
we get

		\begin{equation*}\begin{split}\label{gaffneyHH h=n}
\| \psi^\sharp&\alpha \|_{W^{2,p}(\tilde B_\eta)}
\le
 C_G\left\{\| 
  \psi^\sharp\alpha\|_{L^p(\he n)}+\| 
   d_c^{\mathbb H}(\psi^\sharp\alpha)\|_{L^p(\he n)} +\| 
   d_c^{\mathbb H}\delta_c^{\mathbb H}(\psi^\sharp\alpha)\|_{L^p(\he n)}\right\}
	\\&
	\le
 C_G\left\{\| 
  \psi^\sharp\alpha\|_{L^p(\tilde B_\eta)}+\| 
  \psi^\sharp (d_c^{M}\alpha)\|_{L^p(\tilde B_\eta)} +\| 
  \psi^\sharp (d_c^{M}\delta_c^M\alpha)\|_{L^p(\tilde B_\eta)}+\| 
   \mathcal Q(\psi^\sharp\alpha)\|_{L^p(\tilde B_\eta)}\right\}
\\&
\le
 C_G\left\{\| 
  \psi^\sharp\alpha\|_{L^p(\tilde B_\eta)}+\| 
  \psi^\sharp (d_c^{M}\alpha)\|_{L^p(\tilde B_\eta)} +\| 
  \psi^\sharp (d_c^{M}\delta_c^M\alpha)\|_{L^p(\tilde B_\eta)}+\varepsilon\| 
   \psi^\sharp\alpha\|_{W^{2,p}(\tilde B_\eta)}\right\}
\end{split}\end{equation*}
Choosing $\varepsilon<1/(2C_G)$ we can absorb the term $\varepsilon\| 
   \psi^\sharp\alpha\|_{W^{2,p}(\tilde B_\eta)}$  in the left hand side and eventually get  \eqref{local ineq to prove due}. 
	
	\medskip
	
	The case $h=n+1$ can be handled similarly, taking into account \eqref{T} and proving therefore \eqref{local ineq to prove tre}.
	
\end{proof}
 As noticed in Remark \ref{covering ad hoc},  if $\alpha$ is supported in $B(a,\eta)$ then   the norms
$$
\| \alpha\|_{W^{\ell,p}(M,E_0^\bullet)} \qquad \mbox{and}\qquad \|\psi^\sharp \alpha\|_{W^{\ell,p}(\he n,E_0^\bullet)}
$$
are equivalent, with equivalence constants independent of $\psi$. From the previous theorem we immediatly get the following local result in $M$.

\begin{remark}\label{lemma ivaniec_M}Using the same notation as the previous theorem,
let $B(a,\eta)$ be a ball satisfying Remark \ref{covering ad hoc}. Let $\alpha$ be a smooth form supported in $B(a,\eta)$. If  $\eta<\tilde\eta(C_M, C_G)$ then there exists a constant $C>0$  depending only on $C_M$  and $C_G$, so that, if $h\neq n, n+1$
\begin{align}
			\Vert\alpha\Vert_{W^{1,p}(B(a,\eta), E_0^{h})}\le C\big(\Vert\alpha\Vert_{L^p(B(a,\eta), E_0^{h})}+\Vert d_c^M\alpha\Vert_{L^p(B(a,\eta), E_0^{h+1})}+\Vert\delta_c^M\alpha\Vert_{L^p(B(a,\eta), E_0^{h-1})}\big)\,.
		\end{align}
		Whereas, if $h= n$ we get
		\begin{align}
			\Vert\alpha\Vert_{W^{2,p}(B(a,\eta),E_0^{n})}\le C\big(\Vert\alpha\Vert_{L^p(B(a,\eta), E_0^n)}+\Vert d_c^M\alpha\Vert_{L^p(B(a,\eta), E_0^{n+1})}+\Vert d_c^M\delta_c^M\alpha\Vert_{L^p(B(a,\eta),E_0^{n})}\big)\,,
		\end{align}
		and if  $h= n+1$ we get
		\begin{align}
			\Vert\alpha\Vert_{W^{2,p}(B(a,\eta), E_0^{n+1})}\le C\big(\Vert\alpha\Vert_{L^p(B(a,\eta),E_0^{n+1})}+\Vert \delta_c^Md_c^M\alpha\Vert_{L^p(B(a,\eta), E_0^{n+1})}+\Vert\delta_c^M\alpha\Vert_{L^p(B(a,\eta), E_0^{n})}\big)\,.
		\end{align}
\end{remark}

\bigskip

We are now in position to prove the following  Sobolev-Gaffney type inequalities on $M$ if we assume $M$ to be a smooth sub-Riemannian contact manifold  without boundary with bounded geometry. 
\begin{theorem}\label{gaffney-M} Let $(M,H,g)$ be a smooth contact manifold   with bounded geometry,  without boundary. Let $1\le h\le 2n$,  and $1<p<\infty$. There
exists a positive constant $C=C(C_M, C_G)$ such that for all $\alpha \in  W^{1,p}(M, E_0^h)$ we have: 
\begin{itemize}

\item[i)] for  $h\neq n, n+1$,
\begin{equation}\label{p>1-M}\begin{split}
\| \alpha \|_{W^{1,p}(M, E_0^h)} &
\le C\big( \| d_c^M\alpha \|_{L^{p}(M, E_0^{h+1})} +
\| \delta_c^M \alpha \|_{L^{p}(M, E_0^{h-1})}
\\& \hphantom{xxxxx}+ \|\alpha\|_{L^{p}(M, E_0^h)}
\big);
\end{split}\end{equation}
\end{itemize}


\begin{itemize}
\item[ii)]  for $h= n$,  
\begin{equation}\label{treno1}\begin{split}
\| \alpha \|_{W^{2,p}(M, E_0^n)} &
\le C\big( \| d_c^M\alpha \|_{L^{p}(M, E_0^{n+1})} +
\| d_c^M\delta_c^M \alpha \|_{L^{p}(M, E_0^{n})}
\\&
 \hphantom{xxxxx} +  \|\alpha\|_{L^p(M, E_0^n)}\big);
\end{split}\end{equation}

\item[iii)]  for $h=  n+1$,   
\begin{equation}\label{p>1:eq2M}\begin{split}
\| \alpha \|_{W^{2,p}(M, E_0^{n+1})} &
\le C\big( \| \delta_c^Md_c^M\alpha \|_{L^{p}(M, E_0^{n+1})} +
\| \delta_c^M \alpha \|_{L^{p}(M, E_0^{n})}
\\&
 \hphantom{xxxxx} +  \|\alpha\|_{L^p(M, E_0^{n+1})}\big).
\end{split}\end{equation}
\end{itemize}


\end{theorem}

		\begin{proof} Let $\frac12{\tilde\eta(C_M, C_G)}<\eta<\tilde\eta(C_M, C_G)$ (where $\tilde\eta$ as in Theorem \ref{lemma ivaniec}) and consider the countable, locally finite, atlas $\mc U:=\{ B(a_j, \eta), \psi_{j}\}$ of Remark \ref{covering ad hoc}, were $\psi_j: B(e,\eta/C_M)\to M$. 
 As in Definition \ref{carte}, 
let now $\{\chi_j\}$ be a partition of the unity {  subordinate} to the atlas.
Without loss of
generality, we can assume $\psi_j^{-1}(\mathrm{supp}\;\chi_j) \subset B(e,\eta/C_M)$.

We have
$$
\alpha = \sum_j \chi_j \alpha
$$
Notice that  $\chi_j\alpha$ is supported in $\psi_j(B(e,\eta/C_M))$. Remember, by definition, for $\ell=1,2$ 
$$
\|\alpha\|_{W^{\ell,p}(M, E_0^\bullet)}:= \left(\sum_j \|\psi_{j}^\sharp (\chi_j\alpha)\|^p_{W^{\ell,p}(\he n, E_0^\bullet)}\right)^{1/p}.
$$
\medskip

In the sequel $c$ will denote a geometric  constant that may vary from line to line, depending  in principle on $C_M$, $C_G$, $\eta$  (and on $p, \ h, \ n$). Once we have chosen $\frac12{\tilde\eta(C_M, C_G)}<\eta<\tilde\eta(C_M, C_G)$, the dependence of $c$ is only on $C_M$, $C_G$ (and on $p, \ h, \ n$).
%


\bigskip

\noindent
$\bullet$ Suppose first $h\neq n,n+1$. We divide the proof in $3$ steps.

\begin{itemize}
	\item[Step 1.] 

Let $j\in \N$ be fixed, and let ($B(a_j,\eta)$, $\psi_j$) be a chart of $\mc U$.
%
 %
%
We apply  Theorem \ref{lemma ivaniec} to $\psi_j^\sharp(\chi_j \alpha)$ 
(see\eqref{local ineq to prove}).

Hence, 
\begin{equation*}
\begin{split}
\|\psi_j^\sharp(\chi_j \alpha)&\|_{W^{1,p}(B(e,\frac{\eta}{C_M}), E_0^h)}
\le c\left\{\|\psi_j^\sharp(d_c^M\chi_j \alpha)\|_{L^p(B(e,\frac{\eta}{C_M}), E_0^{h+1})} \right.
\\
&\left.+\| 
   \psi_j^\sharp(\delta_c^M\chi_j \alpha)\|_{L^p(B(e,\frac{\eta}{C_M}), E_0^{h-1})}+\| 
  \psi_j^\sharp(\chi_j \alpha)\|_{L^p(B(e,\frac{\eta}{C_M}), E_0^{h})}\right\}\,.
	\end{split}
\end{equation*}

Now, since
$$d_c^M(\chi_j\alpha)= \chi_jd_c^M\alpha+[d_c^M,\chi_j]\alpha, 
\quad\quad \delta_c^M(\chi_j\alpha)= \chi_j\delta_c^M\alpha+[\delta_c^M,\chi_j]\alpha\, ,$$ 
 from the previous inequality we get

\begin{align}\label{fney_Mm}
\begin{split}
\|\psi_j^\sharp(\chi_j \alpha)&\|^p_{W^{1,p}(B(e,\frac{\eta}{C_M}), E_0^h)}\le
	c\left\{\| \psi_j^\sharp(
   \chi_jd_c^M\alpha)\|^p_{L^p(B(e,\frac{\eta}{C_M}), E_0^{h+1})} \right.
	\\&
	\left.+\| 
   \psi_j^\sharp(\chi_j\delta_c^M\alpha)\|^p_{L^p(B(e,\frac{\eta}{C_M}), E_0^{h-1})}+\| 
  \psi_j^\sharp(\chi_j \alpha)\|^p_{L^p(B(e,\frac{\eta}{C_M}), E_0^{h})}\right.
	\\&
	\left. + \| \psi_j^\sharp([d_c^M,\chi_j]\alpha)\|^p_{L^p(B(e,\frac{\eta}{C_M}), E_0^{h+1})}+\| \psi_j^\sharp([\delta_c^M,\chi_j]\alpha)\|^p_{L^p(B(e,\frac{\eta}{C_M}), E_0^{h-1})}
	\right\}
	\end{split}
\end{align}



\item[Step 2.]

We show now that we can control the sum with respect to $j$, of the last two terms in \eqref{fney_Mm}, with the norm $\|\alpha\|^p_{L^p(M,E_0^h)} $.

First, by Lemma \ref{commutator bis}, \eqref{pullback commutator i.} and \eqref{pullback commutator ii.},

$$ 
   \psi_j^\sharp[d_c^M,\chi_j]\alpha=[d_c^\mathbb H, \chi_j\circ \psi_j]\psi_j^\sharp\alpha ,$$ and 
	
	$$ 
   \psi_j^\sharp[\delta_c^M,\chi_j]\alpha=[\delta_c^\mathbb H, \chi_j\circ \psi_j]\psi_j^\sharp\alpha+\ast [d_c^\mathbb H, \chi_j\circ \psi_j]\psi_j^\sharp\ast\alpha.$$

On the other hand, by Lemma \ref{leibniz}, the 
differential operators   $[d_c^\mathbb H,\chi_j\circ \psi_j]$ and $[\delta_c^\mathbb H,\chi_j\circ \psi_j]$ in $\he n$ 
have order  0 if $h \neq n, n+1$. 

	Keeping in mind this fact, we start from the estimate of the $L^p$- norm of $\psi_j^\sharp[d_c^M,\chi_j]\alpha$. We have:
	
	\begin{equation}
	\begin{split}	
	\|[d_c^\mathbb H, \chi_j\circ \psi_j]\psi_j^\sharp\alpha\|_{L^p(B(e,\frac{\eta}{C_M}), E_0^{h+1})}&\le c \sum_{k\in I_j}\| [d_c^\mathbb H, \chi_j\circ \psi_j]\psi_j^\sharp(\chi_k\alpha)\|_{L^p(B(e,\frac{\eta}{C_M}), E_0^{h+1})}
	\\&
	\le c  \sum_{k\in I_j}\|\psi_j^\sharp(\chi_k\alpha)\|_{L^p(B(e,\frac{\eta}{C_M}), E_0^{h})}
	\\&
	\le c  \sum_{k\in I_j}\|\psi_j^\sharp(\psi_k^\sharp)^{-1}\psi_k^\sharp(\chi_k\alpha)\|_{L^p(B(e,\frac{\eta}{C_M}), E_0^{h})}
	\\&
	\le c  \sum_{k\in I_j}\|\psi_k^\sharp(\chi_k\alpha)\|_{L^p(B(e,\frac{\eta}{C_M}), E_0^{h})},
	\end{split}	
	\end{equation}
	where $\#I_j\le N(\eta)$ since, by Remark \ref{covering ad hoc}, $B(a_j,\eta)$ intersect  at most $N(\eta)$ other balls of
the covering (and where the constant also depends on the uniform bound by $C_M$, of the horizontal derivatives of $\psi_k^{-1}\circ\psi_j$). However, since $\eta$ is bounded  above and below by a quantity that depends only on $C_M$ and $C_G$, also the $\#I_j$ can be controlled only by a geometric constant. 
Hence, also
\begin{equation*}
	\begin{split}	
	\|[d_c^\mathbb H, \chi_j\circ \psi_j]\psi_j^\sharp\alpha\|^p_{L^p(B(e,\frac{\eta}{C_M}), E_0^{h+1})}&
	\le c  \sum_{k\in I_j}\|\psi_k^\sharp(\chi_k\alpha)\|^p_{L^p(B(e,\frac{\eta}{C_M}), E_0^{h})}.
	\end{split}	
	\end{equation*}
	Finally, an analogous estimate of the $L^p$- norm of $[\delta_c^\mathbb H, \chi_j\circ \psi_j]\psi_j^\sharp\alpha$, gives an estimate of the $L^p$- norm of $\psi_j^\sharp[\delta_c^M,\chi_j]\alpha$ . Indeed,
	\begin{equation*}
	\begin{split}	
	\|[\delta_c^\mathbb H, \chi_j\circ \psi_j]\psi_j^\sharp\alpha\|^p_{L^p(B(e,\frac{\eta}{C_M}), E_0^{h-1})}&
	\le c  \sum_{k\in I_j}\|\psi_k^\sharp(\chi_k\alpha)\|^p_{L^p(B(e,\frac{\eta}{C_M}), E_0^{h})}.
	\end{split}	
	\end{equation*}

	Using, once again, the fact that the cover is uniformly locally finite and $\sum_j\sum_{k\in I_j}=\sum_{k}\sum_{j\in I_k}$, summing up over $j$ in the inequalities above and keeping in mind Remark \ref{norme_lp}, we get 
	 
	\begin{equation*}
	\begin{split}
	\sum_j\| \psi_j^\sharp&([d_c^M,\chi_j]\alpha)\|^p_{L^p(B(e,\frac{\eta}{C_M}), E_0^{h+1})}+\sum_j\| \psi_j^\sharp([\delta_c^M,\chi_j]\alpha)\|^p_{L^p(B(e,\frac{\eta}{C_M}), E_0^{h-1})}
	\\&\le c \sum_k\| \psi_k^\sharp(\chi_k
    \alpha)\|^p_{L^p(B(e,\frac{\eta}{C_M}), E_0^{h})}\le c \| 
    \alpha\|^p_{L^p(M, E_0^{h})}
	\end{split}	
	\end{equation*}

	\item[Step 3.]	Finally, summing up over $j$ in  \eqref{fney_Mm},  and using the last estimates, we get 
\begin{align*}
\begin{split}
	&\left(\sum_j\|\psi_j^\sharp(\chi_j\alpha)\|^p_{W^{1,p}(B(e,\frac{\eta}{C_M}), E_0^h)}\right)^{1/p}
\\&
\le c  \left\{\left(
    \sum_j \|\psi_j^\sharp(\chi_j d_c^M\alpha )\|^p_{L^p(B(e,\frac{\eta}{C_M}), E_0^{h+1})} \right)^{1/p}+\left(
    \sum_j\|\psi_j^\sharp(\chi_j \delta_c^M\alpha)\|^p_{L^p(B(e,\frac{\eta}{C_M}), E_0^{h-1})} \right)^{1/p}
  \right.
\\&
+\left.
    \left(
    \sum_j\| 
  \psi_j^\sharp(\chi_j\alpha)\|^p_{L^p(B(e,\frac{\eta}{C_M}), E_0^h)}\right)^{1/p}\right\}+\|\alpha\|_{L^p(M, E_0^{h})}\,.
\end{split}
\end{align*}
which, keeping again in mind Remark \ref{norme_lp}, gives eventually \eqref{p>1-M}.

\end{itemize}

 
		\bigskip
		
	\noindent $\bullet$	Let now $h=n$. 
	
	The argument above needs to be  slightly modified. For $j\in\N$ fixed and keeping in mind again that $\chi_j\alpha$ is supported in $\psi_j(B(e,\frac{\eta}{C_M}))$, we  apply now \eqref{local ineq to prove due} to get

		\begin{equation}
		\begin{split}\label{gaffneyM-h=n}
&
 \|\psi_j^\sharp (\chi_j\alpha)\|^p_{W^{2,p}(B(e,\frac{\eta}{C_M}),  E_0^n)}\le c\left\{\| 
   \psi_j^\sharp(\chi_j d_c^M\alpha)\|^p_{L^p(B(e,\frac{\eta}{C_M}),  E_0^{n+1})} +\right.
\\&
\left.\| 
   \psi_j^\sharp(\chi_jd_c^M\delta_c^M \alpha)\|^p_{L^p(B(e,\frac{\eta}{C_M}),  E_0^{n})}+\| 
  \psi_j^\sharp(\chi_j\alpha)\|^p_{L^p(B(e,\frac{\eta}{C_M}),  E_0^{n})}
  \right.
\\&
+\left. 
\| 
   \psi_j^\sharp([d_c^M,\chi_j]\alpha)\|^p_{L^p(B(e,\frac{\eta}{C_M}),  E_0^{n+1})} +
\| 
   \psi_j^\sharp([d_c^M\delta_c^M,\chi_j]\alpha)\|^p_{L^p(B(e,\frac{\eta}{C_M}),  E_0^{n})} \right\}\,.
\end{split}\end{equation}

		Since $h=n$, by Lemma \ref{leibniz} and Lemma \ref{commutator bis}, $\psi_j^\sharp[d_c^M,\chi_j]$ and $\psi_j^\sharp[d_c^M\delta_c^M,\chi_j]$  are now operators of order $1$. Reasoning as in Step 2, we can write $$\| 
  \psi_j^\sharp([d_c^M,\chi_j]\alpha)\|^p_{L^p(B(e,\frac{\eta}{C_M}), E_0^{n+1})}\le \sum_{k\in I_j}\| 
  \psi_k^\sharp(\chi_k\alpha)\|^p_{W^{1,p}(B(e,\frac{\eta}{C_M}), E_0^{n})}\,.$$ If we combine this with Remark \ref{magenes}, then for any $0<\varepsilon<1$ there exists a constant $c(\varepsilon)$ such that   
\begin{equation}\label{quasi}
\begin{split}
	\sum_j&\| 
  \psi_j^\sharp([d_c^M,\chi_j]\alpha)\|^p_{L^p(B(e,\frac{\eta}{C_M}), E_0^{n+1})} 
	\le \sum_k\|\psi_k^\sharp(\chi_k\alpha)\|^p_{W^{1,p}(B(e,\frac{\eta}{C_M}),  E_0^{n})}
	\\&
		\le \varepsilon \sum_k\| 
    \psi_k^\sharp   (\chi_k\alpha)\|^p_{W^{2,p}(B(e,\frac{\eta}{C_M}),  E_0^{n})} +c(\varepsilon)\|
    \sum_k\phi_k^\sharp  (\chi_k\alpha)\|^p_{L^p(B(e,\frac{\eta}{C_M}),  E_0^{n})}
		\\&
	\le \varepsilon 	\|\alpha\|^p_{W^{2,p}(M,  E_0^{n})} +c(\varepsilon)\|\alpha\|^p_{L^{p}(M,  E_0^{n})}\,.
\end{split}
\end{equation}
		A similar argument shows that
\begin{equation}\label{quasiquasi}
\begin{split}\sum_j&\| 
  \psi_j^\sharp( [d_c^M\delta_c^M,\chi_j](\chi_j\alpha))\|^p_{L^p(B(e,\frac{\eta}{C_M}),E_0^n)}
	\le \varepsilon 	\|\alpha\|^p_{W^{2,p}(M,  E_0^{n})} +c(\varepsilon)\|\alpha\|^p_{L^{p}(M,  E_0^{n})}\,.
		\end{split}
\end{equation}
		
	Going back to \eqref{gaffneyM-h=n}, summing over $j$ and taking the power $1/p$ of all the addends, we get
	\begin{equation*}\label{treno11}
	\begin{split}
\| \alpha \|_{W^{2,p}(M, E_0^n)} &
\le c\big( \| d_c^M\alpha \|_{L^{p}(M, E_0^{n+1})} +
\| d_c^M\delta_c^M \alpha\|_{L^{p}(M, E_0^{n})}
\\&
 +  \|\alpha\|_{L^p(M, E_0^n)}+ \varepsilon 	\|\alpha\|_{W^{2,p}(M,  E_0^{n})} +c(\varepsilon)\|\alpha\|_{L^{p}(M,  E_0^{n})}\big).
\end{split}
\end{equation*}

	Therefore, absorbing 	$\varepsilon	\|\alpha\|_{W^{2,p}(M,  E_0^{n})}$ in the left hand side, up to changing the constants from line to line we  get \eqref{treno1}.
	
	\bigskip
		
\noindent $\bullet$		The case $h=n+1$.

 We fix again $j\in\N$ and we start from the local estimate \eqref{local ineq to prove tre}. The case $h=n+1$ can be dealt with a similar argument to the case $h=n$. Indeed, we have to use the fact that the differential operators $\psi_j^\sharp[\delta^M_cd_c^M,\chi_j]$ and $\psi_j^\sharp[\delta_c^M,\chi_j]$, again by Lemma \ref{commutator bis} and Lemma \ref{leibniz},  are  operators in $\he n$ of order $1$. Then we can produce estimates for this operators analogous to the ones in \eqref{quasi} and \eqref{quasiquasi}. Finally, to conclude we need, once again, Remark \ref{magenes}. 
		\end{proof}

		\bibliographystyle{amsplain}

\bibliography{BTTr_submitted}

\bigskip
\tiny{
\noindent
Annalisa Baldi and Maria Carla Tesi 
\par\noindent
Universit\`a di Bologna, Dipartimento
di Matematica\par\noindent Piazza di
Porta S.~Donato 5, 40126 Bologna, Italy.
\par\noindent
e-mail:
annalisa.baldi2@unibo.it, 
mariacarla.tesi@unibo.it.
}

\medskip

\tiny{
\noindent
Francesca Tripaldi 
\par\noindent Mathematisches Institut,
\par\noindent University of Bern, 
\par\noindent Sidlerstrasse 5
3012 Bern
, Switzerland.
\par\noindent
e-mail: francesca.tripaldi@unibe.ch
}

\end{document}